\documentclass{amsart}
\usepackage[all]{xy}
\usepackage{amscd}
\numberwithin{equation}{section}
\newtheorem*{atiyahexample}{Example}
\let\cal\mathcal
\def\Ascr{{\cal A}}

\def\Cscr{{\cal C}}
\def\Dscr{{\cal D}}
\def\Escr{{\cal E}}
\def\Fscr{{\cal F}}

\def\Hscr{{\cal H}}

\def\Kscr{{\cal K}}
\def\Lscr{{\cal L}}
\def\Mscr{{\cal M}}
\def\Nscr{{\cal N}}
\def\Oscr{{\cal O}}
\def\Pscr{{\cal P}}
\def\Qscr{{\cal Q}}
\def\Rscr{{\cal R}}
\def\Sscr{{\cal S}}
\def\Tscr{{\cal T}}

\let\blb\mathbb

\def \PP{{\blb P}}

\def \ZZ{{\blb Z}}

\def \NN{{\blb N}}

\def\Id{\operatorname{id}}

\def\Lotimes{\overset{L}{\otimes}}

\def\mod{\operatorname{mod}}

\def\Qch{\mathop{\text{Qch}}}
\def\coh{\mathop{\text{\upshape{coh}}}}

\def\rad{\operatorname {rad}}

\def\Spec{\operatorname {Spec}}

\def\depth{\operatorname {depth}}

\def\Ext{\operatorname {Ext}}
\def\Hom{\operatorname {Hom}}
\def\End{\operatorname {End}}
\def\RHom{\operatorname {RHom}}
\def\uRHom{\operatorname {R\mathcal{H}\mathit{om}}}
\def\Sl{\operatorname {Sl}}

\def\Tr{\operatorname {Tr}}
\def\coker{\operatorname {coker}}
\def\ker{\operatorname {ker}}

\def\End{\operatorname {End}}

\def\rk{\operatorname {rk}}

\def\Pic{\operatorname {Pic}}

\def\r{\rightarrow}

\DeclareMathOperator{\Proj}{Proj}

\newtheorem{lemma}{Lemma}[section]

\newtheorem{theorem}[lemma]{Theorem}

\newtheorem{lemmas}{Lemma}[subsection]
\newtheorem{propositions}[lemmas]{Proposition}
\newtheorem{theorems}[lemmas]{Theorem}
\newtheorem{corollarys}[lemmas]{Corollary}

\theoremstyle{definition}

\newtheorem{example}[lemma]{Example}

{

}

\theoremstyle{remark}

\newtheorem{remark}[lemma]{Remark}
\newtheorem{remarks}[lemmas]{Remark}

\newdimen\uboxsep \uboxsep=1ex
\def\uboxn#1{\vtop to 0pt{\hrule height 0pt depth 0pt\vskip\uboxsep
\hbox to 0pt{\hss #1\hss}\vss}}

\def\uboxs#1{\vbox to 0pt{\vss\hbox to 0pt{\hss #1\hss}
\vskip\uboxsep\hrule height 0pt depth 0pt}}

\def\uHom{\Hscr \mathit{om}}
\def\uEnd{\Escr \mathit{nd}}

\def\red{\operatorname{red}}

\def\Per{\operatorname{Per}}
\def\Cl{\operatorname{Cl}}
\def\Qch{\operatorname{Qch}}

\author{Michel Van den Bergh}
 \address{Departement WNI,  Limburgs Universitair Centrum, 
 3590 Diepenbeek, Belgium.}
  \email{vdbergh@luc.ac.be}
\thanks{The author is a director of research at the FWO}
\date{July 20, 2002}
\title[Three-dimensional flops and non-commutative rings]{Three-dimensional flops and non-commutative rings}
\begin{document}
\begin{abstract}
For $Y,Y^+$ three-dimension\-al  smooth varieties  related by
a flop, 
Bondal and Orlov conjectured  that the derived
categories $D^b(\coh(Y))$ and
 $D^b(\coh(Y^+))$ are equivalent. This conjecture was recently proved by
Bridgeland. Our aim in this paper is to give a partially new
proof of Bridgeland's result using non-commutative rings. The new proof 
also covers some mild singular
and higher dimensional situations (including the one  in the recent
paper by Chen: ``Flops and Equivalences of derived Categories for
Threefolds with only Gorenstein Singularities''). 
\end{abstract}
\maketitle
\section{Introduction}
\label{ref-1-0}
Let $k$ be an algebraically closed field of characteristic zero and
assume that
$Y,Y^+$ are $3$-dimension\-al  smooth varieties over $k$ related by
a flop. 
Bondal and Orlov conjectured in \cite{Bondal1} (and proved in some
special cases) that the derived
categories $D^b(\coh(Y))$ and
 $D^b(\coh(Y^+))$ are equivalent (see also the talk at the IMU
 congress \cite{BO2}). This conjecture was recently proved by
Bridgeland \cite{Br1}. Using similar methods Chen proved that
the derived equivalence also holds in some mild Gorenstein cases \cite{Chen}. Recently
Kawamata used Chen's result to obtain an analogue of the Bondal-Orlov
conjecture which is also valid in the non-Gorenstein case \cite{Kawamata}.

Our aim in this paper is to give a partially new
proof of Bridgeland's result. The new proof 
also covers some mild singular
and higher dimensional situations such as the one  in the
paper by Chen \cite{Chen}. See Theorem \ref{ref-C-5} below.

In the
interest of full disclosure we should mention that 
Bridgeland's methods (and Chen's extension of these) also yield the
existence of flops whereas in our method we have to assume
this. Three-dimensional flops (in contrast to flips!) are very cheap
however (see \cite[Theorem 2.4]{Kollar}).

Saying that $Y,Y^+$ are related by a flop means that there is a diagram 
\begin{equation}
\label{ref-1.1-1}
\xymatrix{
Y\ar[dr]_{f} & &\ar[dl]^{f^+} Y^+\\
& X &
}
\end{equation}
of birational proper maps with fibers of dimension
$\le 1$ such that there is an ample divisor $D$ on $Y$ with the property that if $E$ is
the strict transform of $D$ on $Y^+$ then $-E$ is ample. According to
\cite[Prop. 16.2]{CKM} $X$ has Gorenstein terminal singularities. 

Bridgeland's proof is based on the now familiar technique of
Fourier-Mukai transform.  A very important instance
of this technique is the general approach to the McKay
correspondence introduced in \cite{BKR}. In that paper it is shown
that if $V$ is a 
$3$-dimensional vector space and $G\subset
\Sl(V)$ is a finite group then  there is a derived equivalence
between the skew
group ring $k[V]{*}G$ and   a suitable
$G$-equivariant Hilbert scheme  of $V$ (see also
\cite{KV} for the two-dimensional case).

An essential feature of the McKay correspondence is the appearance of
the skew group algebra $k[V]{*}G$. This suggests that
it should perhaps  be possible to explain Bridgeland's proof in
terms of non-commutative algebra. 
\begin{atiyahexample}
The most trivial flop is the 
so-called ``Atiyah flop''. Let
$R=k[u,v,x,y]/(uv-xy)$ and $X=\Spec R$.  
Let $I$ be the reflexive ideal in $R$ given by $(u,x)$.
The
singular variety $X$ has two resolutions $Y$, $Y^+$ which are obtained by blowing up
either $I$ or 
$I^{-1}\cong (u,y)$. Put
\begin{equation}
\label{ref-1.2-2}
A=\begin{pmatrix} R & I \\
I^{-1} & R
\end{pmatrix}
\end{equation}
It was observed by many people that there are derived equivalences
\[
D^b(\coh(Y))\cong D^b(\mod(A))\cong D^b(\coh(Y^+))
\]
\end{atiyahexample}
So the question is how to construct an analogue of the ring $A$ in
general. Luckily Bridgeland's paper is very helpful: Bridgeland constructs a
series of t-structures on $D^b(\coh(Y))$ indexed by a ``perversity''
$p\in\ZZ$ whose hearts are abelian
categories  denoted by
${}^p\Per(Y/X)$. Below we will always assume $p=0,-1$. In
that case ${}^p\Per(Y/X)$ is simply a ``tilting'' of $\coh(Y)$ in the sense of
\cite{HRS} (see \S\ref{ref-3-6} for precise definitions). We will show that if $X=\Spec R$ is affine then the
category ${}^p\Per(Y/X)$ has a projective generator. Our desired
non-commutative ring is the endomorphism ring of this generator.

In fact this part of the argument works in far greater generality. We
prove the following result (this is a combination of Propositions
\ref{ref-3.3.1-24} and \ref{ref-3.2.7-18})
{\def\thelemma{A}
\begin{theorem} Let $f:Y\r X$ be a projective morphism  between
  quasi-projective schemes over an affine scheme such that the fibers
  of $f$ have dimension $\le 1$ and such that $Rf_\ast
  \Oscr_Y=\Oscr_X$. Then there exists a 
  vector bundle $\Pscr$ 
  on $Y$ with the following properties.
\begin{enumerate}
\item Let $\Ascr=f_\ast \uEnd(\Pscr)$. The functors
$Rf_\ast \uRHom_Y(\Pscr,-)$ and $f^{-1}(-)\Lotimes_{f^{-1}(\Ascr)}
\Pscr$ define inverse equivalences between $D^b(\coh(Y))$ and
$D^b(\coh(\Ascr))$. These equivalences restrict to equivalences
between ${}^{-1}\Per(Y/X)$ and $\coh(\Ascr)$.
\item The functors
$Rf_\ast \uRHom_Y(\Pscr^\ast,-)$ and $f^{-1}(-)\Lotimes_{f^{-1}(\Ascr^\circ)}
\Pscr^\ast$ define inverse equivalences between $D^b(\coh(Y))$ and
$D^b(\coh(\Ascr^\circ))$. These equivalences restrict to equivalences
between ${}^0\Per(Y/X)$ and $\coh(\Ascr^\circ)$.
\end{enumerate}
\end{theorem}
}
The vector bundle $\Pscr$ which occurs in the above theorem is a
so-called \emph{local projective generator} for
${}^{-1}\Per(Y/X)$. I.e. an object in ${}^{-1}\Per(Y/X)$ such that there
exists an affine covering of $X$: $X=\bigcup_i U_i$ with the property
that for all $i$: $\Pscr\mid U_i$  is a projective generator for ${}^{-1}
\Per(f^{-1}(U_i)/U_i)$. 

Now let us return to the situation of  diagram \eqref{ref-1.1-1}. The
hypotheses on the dimension of the fibers yield that the categories of
reflexive sheaves on $Y,Y^+$ and $X$ are equivalent. 
Choose a
vector bundle $\Pscr$ as in Theorem A. Under the
above equivalences we find a corresponding sheaf of reflexive modules
$\Qscr^+$ on $Y^+$. Consulting Bridgeland's paper for inspiration we
note that Bridgeland shows that ${}^{-1}\Per(Y/X)$ is equivalent to
${}^0\Per(Y/X)$. So this suggests
we should
prove that $\Qscr^+$ is a local projective generator for
${}^0\Per(Y^+/X)$.  Indeed if this is the case then by using Theorem A
for $Y$ and $Y^+$ we find
\[
D^b(\coh(Y))\cong D^b(\coh(\Ascr))\cong D^b(\coh(Y^+))
\]
where $\Ascr=f_\ast\uEnd_Y(\Pscr)=f_\ast\uEnd_Y(\Qscr^+)$.

We show that the property of being a local projective generator may be
verified in the completions of the closed points of $X$. 
Hence we may reduce to the case where $X$ is
the spectrum of a complete local ring. In that case we can invoke the
methods of Artin and Verdier \cite{AV} to give a precise
classification of the
projective objects in ${}^p\Per(Y/X)$. Again this works in greater
generality. Our result is the following
(this is  Theorem \ref{ref-3.5.5-45})
{\def\thelemma{B}
\begin{theorem} Assume that we are in the situation of Theorem
  A and assume in addition that $X$ is the spectrum of a complete
  local $k$ algebra $R$ with residue field equal to $k$.  Let $(C_i)_i$ be
  the irreducible components of the special fiber of $f$. Then the map
  $\Lscr\mapsto \deg(\Lscr\mid C_i)_{i=1,\ldots,n}$ defines an
  isomorphism $\Pic(Y)\cong \ZZ^n$.  Let $\Lscr_i\in \Pic(Y)$ be such
  that $\Lscr_i\mid C_j=\delta_{ij}$.  Define $\Mscr_i$ as the
  middle term of the exact sequence
\[
0\r \Oscr_Y^{r-1}\r  \Mscr_i\r \Lscr_i\r 0
\]
associated to a minimal set of generators of $H^1(Y,\Lscr_i^{-1})$ as
$R$-module. Then the indecomposable projective objects in
${}^{-1}\Per(Y/X)$ are precisely the objects
$(\Mscr_i)_{i=1,\ldots,n}$  and $\Oscr_Y$. The projective generators
of ${}^{-1}\Per(Y/X)$ are $\Oscr_Y^{\oplus a_0}\oplus
\bigoplus_{i=1,\ldots,n} \Mscr_i^{\oplus a_i}$ with $a_i>0$ for all $i$. 
Analogous results hold for ${}^0\Per(Y/X)$ if we replace $\Mscr_i$ by
$\Mscr_i^\ast$. 
\end{theorem}
}
We now return to the situation of diagram \eqref{ref-1.1-1} but we assume
that we are in the formal case, i.e. $X=\Spec R$ where $R$ is a
complete local ring.  We will adorn notations with refer to $f^+$ by a
superscript ``$+$''. 

In the formal case there is a trivial description of
$Y^+$ \cite{Kollar}.  Recall that Gorenstein
terminal singularities are 
hypersurface singularities of multiplicity two. Hence the equation of $R$ may
be written as $u^2+\cdots=0$. It follows that  $u\mapsto -u$ defines an
automorphism $\sigma$ of $X$ of order two. We may take $Y^+=Y$
and $f^+=\sigma \circ f$. We put $C^+_i=C_i$.

From Theorem B,  and the discussion preceeding it, we see
that we are done if we can show that $\Mscr_i$ corresponds to 
$(\Mscr_i^{+})^{\ast}$ under the equivalence of  the categories of reflexive
modules on $Y,Y^+$ and $X$.

Let $I_i=\Gamma(Y,\Lscr_i)$. It is easy to see that $I_i$ is a
reflexive  $R$-module of rank one
Put  $M_i=\Gamma(Y,\Mscr_i)$. We show that the
$M_i$ are indecomposable 
Cohen-Macaulay modules on $R$ which occur as middle term of an exact sequence 
\[
0\r R^{r_i-1} \r M_i \r I_i \r 0
\]
associated to a minimal set of generators of $\Ext^1_R(I_i,R)$ as
$R$-module. In addition we also show that there are exact sequences 
\begin{equation}
\label{ref-1.3-3}
0\r M_i^{\ast} \r  R^{r_i+1} \r I_i \r 0
\end{equation}

Let $\Cl(R)$ be the divisor class group of $R$. The automorphism
$\sigma$ induces the operation $I\mapsto I^{-1}$ on $\Cl(R)$ and from
this we deduce $\Gamma(Y^+,\Lscr_i^+)=I_i^{-1}$. Thus $M_i^+$ occurs 
as middle term of an exact sequence 
\begin{equation}
\label{ref-1.4-4}
0\r R^{r^+_i-1} \r M^+_i \r I^{-1}_i \r 0
\end{equation}
associated to a minimal set of generators of $\Ext^1_R(I^{-1}_i,R)$ as
$R$-module. 

We show in \S\ref{ref-4.1-47} that the Cohen-Macaulay modules defined by
\eqref{ref-1.3-3} and \eqref{ref-1.4-4} are the
same. This finishes the proof.

As indicated above, our methods are valid in a somewhat more general
setting. The precise result we prove is the following (this is Theorem
\ref{ref-4.4.2-59}) 
{\def\thelemma{C}
\begin{theorem}
\label{ref-C-5}
Let $f:Y\r X$ be a projective birational map between normal
quasi-projective  Gorenstein
$k$-varieties of dimension $n\ge 3$ with fibers of dimension $\le 1$ and
assume that the exceptional 
locus of $f$ has codimension $\ge 2$ in $Y$. Assume that $X$ has canonical
hypersurface singularities of multiplicity $\le 2$.  Let $f^+:Y^+\r X$
be the flop of $f$.  Then
  $D^b(\coh(Y))$ and
  $D^b(\coh(Y^+))$ are equivalent and we may choose  this equivalence 
in such a way that ${}^{-1}\Per(Y/X)$ corresponds to
 ${}^{0}\Per(Y^+/X)$.
\end{theorem}
}
To finish this introduction we make some remarks of a more
philosophical nature. If we return to the McKay correspondence then we
see that  the singular variety $V/G=\Spec k[V]^G$ has
two crepant ``resolutions'', a commutative one given by the
$G$-equivariant Hilbert scheme and a
non-commutative one given by the skew group ring $k[V]{*}G$.  Both
resolutions are derived equivalent, so in some sense it doesn't
matter which one we take. 

If $\dim V\ge 4$ then there are examples where $V/G$ does not have a
commutative crepant resolution. However the non-commutative resolution given
by $k[V]{*}G$ of course always exists. So it seems that at least in some
situations non-commutative resolutions are strictly more general than
commutative ones. 

An obvious question is whether something similar is true for a variety
$X$ 
with three-dimensional terminal Gorenstein singularities.  I.e. does
there always exists a ``crepant'' non-commutative resolution of
$X$? In Appendix \ref{ref-A-61} I give a counter example that shows that
this is not the case. I actually think that for three-dimensional terminal
Gorenstein singularities the existence of commutative and
non-commutative crepant resolutions are equivalent. This can presumably be
proved with the same Fourier-Mukai method which was used to establish
the three-di\-men\-si\-on\-al McKay correspondence. 

Most of the results in this paper were conceived during the
2002 OberWolfach meeting suitably called ``Interactions between Algebraic
Geometry and Noncommutative Algebra''. During that meeting, and also at
other times, the author has greatly benefited from discussions with
Alexei Bondal, Tom Bridgeland, Alastair King and Aidan Schofield.
\section{Notations and conventions}
$k$ will always be an algebraically closed field of characteristic
zero. Schemes are \emph{not} supposed to be $k$-schemes, unless this is
explicitly stated. The characteristic zero hypotheses on $k$ could be
avoided at the cost of more technicality in the statement of some
results (i.e. those related to rational singularities). 

If $A$ is a noetherian ring then $\mod(A)$ is the category of finitely
generated right $A$-modules and $D^b(A)=D^b(\mod(A))$. 
If $X$ is a noetherian scheme then $\Qch(X)$ and $\coh(X)$ are the
categories  of quasi-coherent and coherent $\Oscr_X$-modules. If
$\Ascr$ is a sheaf of $\Oscr_X$ algebras then $\coh(\Ascr)$ is the
category of right coherent $\Ascr$-modules.

If $X$ is a scheme and $D$ is a closed subscheme of $X$ then
$\Oscr_X(-D)$ denotes the ideal sheaf of $D$. If $\Oscr_X(-D)$ is
invertible then $\Oscr_X(nD)=\Oscr_X(-D)^{\otimes -n}$.

All Cohen-Macaulay modules in this paper are maximal Cohen-Macaulay.

\section{Acyclic morphisms with one dimensional fibers}
\label{ref-3-6}
\subsection{Generalities}
\label{ref-3.1-7}
Let  $f:Y\r X$ be a projective map between noetherian schemes 
We impose  the following conditions:
\begin{enumerate}
\item $Rf_\ast \Oscr_Y=\Oscr_X$.
\item  The fibers  of $f$ are one-dimensional.
\end{enumerate}
At this stage we do not assume that $f$ is birational.  Note
however that by \cite[Cor. II.11.3]{H} the fibers of $f$ are
connected. 

We consider certain categories of perverse coherent sheaves 
introduced by Bridgeland in \cite{Br1}. 
 We let $\Cscr$ be the abelian(!) subcategory of $\coh(Y)$
consisting  of objects $E$ such that $Rf_\ast E=0$ (this is a
deviation of the notation used by Bridgeland). The following lemma was
proved by Bridgeland.
\begin{lemmas} \label{ref-3.1.1-8} For $F\in D^b(\coh(Y))$ one has $Rf_\ast F=0$ if and
  only if $H^i(F)\in \Cscr$ for all $i$.
\end{lemmas}

We define the following torsion theories on 
$\coh(Y)$. 
\begin{align*}
\Tscr_{-1}&=\{T\in \coh(Y)\mid R^1f_\ast(T)=0, \Hom(T,\Cscr)=0\}\\
\Fscr_{-1}&=\{F\in \coh(Y)\mid f_\ast(F)=0\}
\end{align*}
\begin{align*}
\Tscr_0&=\{T\in \coh(Y)\mid R^1f_\ast(T)=0 \}\\
\Fscr_0&=\{F\in \coh(Y)\mid f_\ast(F)=0, \Hom(\Cscr,F)=0\}
\end{align*}
Then  on $D=D^b(\coh(Y))$ we consider the associated perverse t-structures.
\begin{align*}
{}^p D^{\le 0}&=\{E\in D^{\le 0}\mid H^0(E)\in \Tscr_p\}\\
{}^p D^{\ge 0}&=\{E\in D^{\ge -1}\mid H^{-1}(E)\in \Fscr_p\}
\end{align*}
for $p=-1,0$. 

We denote the heart of these t-structures by  ${}^p\Per(Y/X)$. I.e ${}^p\Per(Y/X)$  consists of the objects $E$ in $D^b(\coh(Y))$
whose only homology lies in degree $-1,0$ such that $H^{-1}(E)\in
\Fscr_p$ and $H^0(E)\in \Tscr_p$. 

It is not entirely clear that  the categories $\Tscr_{-1}$
and $\Fscr_0$ are compatible with restriction since their definition involves
the functor
$\Hom$ instead of $\uHom$. This defect is taken care of by the next
two lemmas.
\begin{lemmas} 
\label{ref-3.1.2-9}
The objects in $\Tscr_{-1}$ are precisely the coherent sheaves $T$ 
  such that the map $f^\ast f_\ast T\r T$ is surjective.
\end{lemmas}
\begin{proof} We concentrate on the non-obvious direction. 
Let $T\in \coh(Y)$ be such that $R^1f_\ast T=0$ and let $T_0$ be the image of $f^\ast f_\ast T\r
 T$.
 We claim $T/T_0\in \Cscr$.  Hence if $T\in\Tscr_{-1}$ then $T=T_0$ and
 thus the canonical map $f_\ast f^\ast T\r T$ is surjective.

We prove the claim. It is clear
that $R^1f_\ast(T/T_0)=0$ and that $f_\ast(T)\r f_\ast(T/T_0)$ is
surjective. So if $T/T_0\not\in \Cscr$ then $f_\ast (T/T_0)\neq
 0$. This means $f^\ast  f_\ast (T/T_0)\r T/T_0$ is not the zero map.
 Hence the composition $\phi:f^\ast f_\ast T \r f^\ast  f_\ast
 (T/T_0)\r T/T_0$ is also not the zero map (since the first map is
 surjective). Now $\phi$ is also the composition of $f^\ast f_\ast T\r
 T \r T/T_0$ which was zero by hypotheses. This contradiction finishes
 the proof. 
\end{proof}
 To obtain an analogous statement
for $\Fscr_0$ we recall that 
$Rf_\ast:D(\Qch(Y))\r D(\Qch(X))$ has a right adjoint $f^!$
\cite{Keller1,Neeman1}.  The identity $Rf_\ast Lf^\ast=\Id$ formally
implies $Rf_\ast f^!=\Id$.
The
explicit formulas for $f^!$ in \cite{RD} show that $f^!$ maps
$D(\coh(X))^{\ge 0}$ to $D(\coh(Y))^{\ge -1}$. 

If $E\in \coh(Y)$ then the composition $E\r f^! Rf_\ast E\r f^!
((R^1f_\ast E)[-1])$ yields a canonical map $\phi_E:E\r H^{-1}(f^!
R^1f_\ast E)$. 
\begin{lemmas} 
The objects in $\Fscr_0$ are precisely the coherent sheaves $F$ 
  such that $\phi_F$ is injective.
\end{lemmas}
\begin{proof} Assume first that $\phi_F$ is
  injective.
If $E\in \coh(X)$ then the fact that $Rf_\ast f^! E=E$ implies $f_\ast
H^{-1}(f^! E)=0$. Hence if $\phi_F$ is injective then $f_\ast
F=0$. Thus $\phi_F$ is the map $F\r H^{0}(f^! Rf_\ast F)$ 
obtained by applying $H^0$ to the map $F\r f^! Rf_\ast F$
given by
adjointness. Assume that there is a non-zero homomorphism $C\r F$ for
$C\in \Cscr$. This yields a non-zero homomorphism $C\r H^{0}(f^!
Rf_\ast F)$ and hence an non-zero homomorphism $C\r f^!
Rf_\ast F$, but this is impossible by adjointness. 

Now assume that $F\in \Fscr_0$. Since $f_\ast F=0$ the map $\phi_X$ is
again the map $F\r H^{0}(f^! Rf_\ast F)$ 
obtained by applying $H^0$ to the map
given by
adjointness. Let $U$ be the cone of $F\r f^! Rf_\ast F$. We have
$Rf_\ast U=0$ and hence by lemma \ref{ref-3.1.1-8} the homology of $U$
lies in $\Cscr$. Looking at the homology exact sequence we obtain
\[
0\r H^{-1}(U)\r F \xrightarrow{\phi_X} H^{0}(f^! Rf_\ast F)
\]
Since $H^{-1}(U)\in \Cscr$ it follows that $H^{-1}(U)=0$ and hence
that $\phi_X$ is injective.
\end{proof}
We obtain the following consequence
\begin{propositions}  The categories $\Tscr_p$, $\Fscr_p$,
  ${}^pD^{\ge 0}$, ${}^pD^{\le 0}$, ${}^p\Per(Y/X)$ are compatible
  with restriction (and more generally flat base change). Furthermore,
  membership of these categories may be checked locally (even for the
  flat topology).
\end{propositions}
One may also check that  if we put for $U\subset X$ open,
${}^p\Per_{Y/X}(U)={}^p\Per(f^{-1}(U)/U)$ then ${}^p\Per_{Y/X}$ is 
 a stack of abelian categories with exact
restriction functors.

\subsection{The case when the base is affine}
\label{ref-3.2-10}
First we assume that $X=\Spec R$ with $R$ a noetherian ring.  Our aim is
to show that ${}^p\Per(Y/X)$ is a module category. We will use the
following version of Morita theory.
\begin{lemmas}
\label{ref-3.2.1-11}
Let $R$ be a noetherian commutative ring and $\Cscr$ an $R$-linear
category such that for $A,B\in \Cscr$ we have that $\Hom_{\Cscr}(A,B)$ is a
finitely generated $R$-module. Assume that $\Pscr$ is a projective object
in $\Cscr$ such that $\Hom_{\Cscr}(\Pscr,E)=0$ implies $E=0$. Put
$A=\End_{\Cscr}(\Pscr)$. Then the functor $\Hom_{\Cscr}(\Pscr,-)$ defines an
equivalence 
between $\Cscr$ and $\mod(A)$ whose inverse is given by $-\otimes_A \Pscr$
(with obvious notations).
\end{lemmas} 
We will call an object $P$ as in the statement of this lemma a
\emph{projective generator} for $\Cscr$.

We use the following result.
\begin{lemmas} 
\label{ref-3.2.2-12}
Let $f:Y\r X$ be a projective morphism between
  noetherian schemes with fibers of dimension $\le n$. Assume that $X$
  is affine. Let $\Lscr$ be an ample line bundle generated by
  global sections and let
  $a\in \ZZ$.

 If $M$ in $D^b(\coh(Y))$ is such that
 $\Ext^i_Y(\Lscr^{a+j},M)=0$ for all $i$ and for $0\le j\le n$ then $M=0$.
\end{lemmas}
\begin{proof}
This is basically \cite[Lemma 4.2.4]{BondalVdB}. For the convenience
of the reader we repeat the proof in the current setting. Without loss
of generality we may assume $a=-n$. 

Use $\Lscr$ to construct a finite map $Y\r \PP^N_X$. The Koszul complex of 
a polynomial ring in $N+1$ variables leads to a long exact sequence on
$\PP^N_X$: 
\[
0\r \Oscr_{\PP^N_X}(-N-1)\r \cdots \r\Oscr_{\PP^N_X}(-u)^{\bigl(\!\begin{smallmatrix}
N+1\\ u\end{smallmatrix}\!\bigr)}
\r\cdots  \r \Oscr_{\PP^N_X}\r 0
\]
The inverse image on $Y$ yields an exact sequence 
\[
0\r \Lscr^{-N-1}\r \cdots\r (\Lscr^{-u})^{\bigl(\!\begin{smallmatrix}
N+1\\ u\end{smallmatrix}\!\bigr)} \r \cdots \r \Oscr_Y\r 0
\]
Let $U$ be the kernel at  $(\Lscr^{-n-1})^{\bigl(\!\begin{smallmatrix}
N+1\\ n+1\end{smallmatrix}\!\bigr)}$. Then the long exact sequence
\[
0\r U \r  (\Lscr^{-n-1})^{\bigl(\!\begin{smallmatrix}
N+1\\ n+1\end{smallmatrix}\!\bigr)}\r \cdots\r (\Lscr^{-1})^{N+1} \r \Oscr_Y\r 0
\]
represents an element of $\Ext^{n+1}_Y(\Oscr_Y,U)$ which must be zero. It
follows that $\Oscr_Y$ is a direct summand in $D^b(\coh(Y))$ of the
complex 
\[
( \Lscr^{-n-1})^{\bigl(\!\begin{smallmatrix}
N+1\\ n+1\end{smallmatrix}\!\bigr)}\r \cdots \r (\Lscr^{-1})^{N+1}
\]
Dualizing and tensoring with  $\Lscr^{-n-p-1}$ for $p\ge 0$ we deduce that
$\Lscr^{-n-p-1}$ is a 
direct summand  of 
\[
(\Lscr^{-n-p})^{N+1} \r \cdots\r ( \Lscr^{-p})^{\bigl(\!\begin{smallmatrix}
N+1\\ n+1\end{smallmatrix}\!\bigr)}
\]
By induction on $p$ we now easily deduce from the hypotheses that
$\Ext^i_Y(\Lscr^j,M)=$ for all $i$ and $j\le 0$. Applying this with
$i=0$ and using the fact that $\Lscr$ is ample this implies $M=0$.
\end{proof}
Now we revert to our blanket assumptions (but we assume $X=\Spec R$ is
affine). 

Let $\frak{V}$ be the category of  vector bundles $\Mscr$ on $Y$
  generated by global sections such that $H^1(Y,\Mscr^\ast)=0$ and let
  $\frak{V}^\ast=\{\Mscr^\ast\mid \Mscr\in \frak{V}\}$.
\begin{lemmas}
\label{ref-3.2.3-13} If $\Mscr\in \frak{V}$  then   for all
  $E\in {}^{-1}\Per(Y/X)$ and for all $i>0$ we have
  $\Ext^i_Y(\Mscr,E)=0$. In particular $\Mscr$ is a projective object
  in ${}^{-1}\Per(Y/X)$.
\end{lemmas}
\begin{proof}
By lemma \ref{ref-3.1.2-9} we know that $\Mscr\in \Tscr_{-1}$. 
It is also clear that if  $F\in\Fscr_{-1}$ then for $i>0$:
$\Ext^i_Y(\Mscr,F[1])=\Ext^{i+1}_Y(\Mscr,F)=0$. 

 So we need to show
   that $\Ext^i_Y(\Mscr,T)=0$ for $T\in
  \Tscr_{-1}$ and for $i>0$. The case $i>1$ follows from the hypotheses on
   the fibers of $f$ so we assume $i=1$. 
Since by lemma \ref{ref-3.1.2-9} $T$ is generated by global sections, $T$ is a
  quotient of some $\Oscr_Y^l$. Let $\Kscr$ be the kernel. Then
  $\Ext^1_Y(\Mscr,T)=\Ext^2_Y(\Mscr,\Kscr)=0$. Hence we are done.
\end{proof}
\begin{lemmas} 
\label{ref-3.2.4-14} If $\Nscr\in \frak{V}^\ast$
  then  for all
  $E\in {}^0\Per(Y/X)$ and for all $i>0$ we have
  $\Ext^i_Y(\Nscr,E)=0$.
In particular $\Nscr$ is a
  projective object in ${}^0\Per(Y/X)$.
\end{lemmas}
\begin{proof}
The fact that  $H^1(Y,\Nscr)=0$ implies that $\Nscr\in \Tscr_0$. 
It is also clear that  if  $F\in\Fscr_{-1}$ then for $i>0$:
$\Ext^i_Y(\Mscr,F[1])=\Ext^{i+1}_Y(\Mscr,F)=0$. 
So we need to show
   that $\Ext^i_Y(\Nscr,T)=0$ for $T\in
  \Tscr_0$ and for $i>0$. The case $i>1$ follows from the hypotheses
  on the fibers of $f$ so we assume $i=1$. Now $\Ext^1_Y(\Nscr,T)=H^1(Y,\Nscr^\ast \otimes_Y T)$. Since
  $\Nscr^\ast$ is generated by global sections it follows that
  $\Nscr^\ast \otimes_Y T$ is a quotient of a number of copies of
  $T$. Hence $\Nscr^\ast \otimes T$ has vanishing cohomology.
\end{proof}
The following is our main result.
\begin{propositions}
\label{ref-3.2.5-15} There exists a vector bundle $\Pscr\in\frak{V}$ 
 which is a projective generator in
  ${}^{-1}\Per(Y/X)$ 
  and whose dual $\Pscr^\ast$ is a projective generator in
  ${}^0\Per(Y/X)$. 
\end{propositions}
\begin{proof} Let $\Lscr$ be an ample line bundle on $Y$ generated
  by global sections.

Let
  $\Pscr_0$ be given by the extension
\begin{equation}
\label{ref-3.1-16}
0\r \Oscr_Y^{r-1} \r \Pscr_0\r \Lscr\r 0
\end{equation}
associated to a set of generators of $H^1(Y,\Lscr^{-1})$ as $R$-module
and put
$\Pscr=\Pscr_0\oplus \Oscr_Y$. 

Clearly $\Pscr_0\in \frak{V}$ and hence by lemma \ref{ref-3.2.3-13}:
$\Ext^i_Y(\Pscr,-)=0$ for 
$i\neq 0$. 
 Hence we only need to show that $\Pscr$ is a generator.

Assume that $E\in D^b(\coh(Y))$ is such that
$\Ext^i_Y(\Pscr,E)=0$ for all $i$.  Then we deduce from this
that $\Ext^i_Y(\Oscr_X,E)=0$ for all $i$ and $\Ext^i_Y(\Lscr,E)=0$ for
all $i$. The conclusion now follows from lemma \ref{ref-3.2.2-12}.

The proof for $p=0$ is similar.
\end{proof}
Now we give  a more detailed discussion of the projective objects
and the projective generators in ${}^p\Per(Y/X)$. 
\begin{propositions}
\label{ref-3.2.6-17}
The projective objects in ${}^{-1}\Per(Y/X)$ are precisely the objects
in $\frak{V}$. The projective objects in ${}^0\Per(Y/X)$ are precisely
the objects in $\frak{V}^\ast$.
\end{propositions}
\begin{proof} Assume $p=-1$. By lemma \ref{ref-3.2.3-13} we already
  know that the objects in $\frak{V}$ are projective. Since
  $\frak{V}$ is clearly closed under direct sums and direct summands,
  the  converse follows from the fact that $\frak{V}$ contains a
  projective generator for ${}^{-1}\Per(Y/X)$. The case $p=0$ is
  identical. 
\end{proof}
If $\Mscr$ is a vector bundle of rank $r$ on $Y$ then by $c_1(\Mscr)$
we denote the class  of $\wedge^r \Mscr$ in $\Pic(Y)$.
\begin{propositions}
\label{ref-3.2.7-18} 
The projective generators in ${}^{-1}\Per(Y/X)$
  are the objects $\Mscr$ in $\frak{V}$ such that $c_1(\Mscr)$ is ample
  and such
  that  $\Oscr_Y$ is a 
  direct summand of some $\Mscr^{\oplus a}$. 
The projective generators in ${}^0\Per(Y/X)$
  are the objects in $\frak{V}^\ast$ which are dual to projective
  generators in ${}^{-1}\Per(Y/X)$
\end{propositions}
\begin{proof} Assume $p=-1$. We first prove that every projective
  generator is of the form stated. Let $\Mscr$ be a projective
  generator of ${}^{-1}\Per(Y/X)$. Then since $\Oscr_Y$ is a projective
  object in ${}^{-1}\Per(Y/X)$, $\Oscr_Y$ must be a direct summand of
  some $\Mscr^{\oplus a}$. 

Let $\Pscr$ be the projective generator constructed in the proof of
Proposition  \ref{ref-3.2.5-15}. Then there exists a $b\in \NN$ such that
$\Mscr^{\oplus b} =\Pscr\oplus \Pscr'$ with $\Pscr'\in
\frak{V}$. Hence $c_1(\Mscr)^b=c_1(\Pscr)c_1(\Pscr')$. Since
$c_1(\Pscr)$ is ample by construction and $c_1(\Pscr')$ is generated
by global sections we deduce from this that $c_1(\Mscr)$ is ample.

Now we prove the converse. Let $\Mscr$ be as in the statement of the
proposition.  By Proposition \ref{ref-3.2.6-17} we already know that
$\Mscr$ is projective. So we only have to show that $\Mscr$ is a generator.

As in the proof of Propositions \ref{ref-3.2.5-15} we will show
that if $E\in D^b(\coh(Y))=0$ is such that 
if $\Ext^i_Y(\Mscr,E)=0$ for all $i$ then $E=0$. Looking at the
appropriate spectral sequence we see that if $\Ext^i_Y(\Mscr,E)=0$ for
all $i$ then $\Ext^i_Y(\Mscr,H^j(E))=0$ for all $i,j$ (this is similar
to lemma \ref{ref-3.1.1-8}). Thus without loss of generality we may
assume $E\in \coh(Y)$.
To prove that $E=0$ it is
sufficient to prove that for every closed point $x$ in $X$ with
defining ideal $m$ we have $E/mE=0$. Consider $m E$. Since this is
a subobject of $E$ we have $\Hom_Y(\Mscr,mE)=0$. However it is also a quotient
of some $E^{\oplus c}$ so we  have in addition $\Ext^1_Y(\Mscr,mE)=0$. Thus
$\Ext^i_Y(\Mscr,E/mE)=0$ for all $i$. Let $C$ be the fiber of $x$. We have
$0=\Ext^i_Y(\Mscr,E/mE)=\Ext^i_C(\Mscr/m\Mscr,E/mE)$. Hence we are now
reduced to the case where $X$ is the spectrum of a field and without
loss of generality we may assume that this field is algebraically
closed. According to lemma \ref{ref-3.5.1-35} below we now have an exact
sequence
\[
0\r \Oscr_Y^{r-1} \r \Mscr\r \Lscr\r 0
\]
with $\Lscr=c_1(\Mscr)$ ample and generated by global sections. Since
$\Oscr_Y$ is a direct summand of some $\Mscr^{\oplus a}$ we also have
$\Ext^i_Y(\Oscr_Y,E)=0$ for all $i$. We now conclude as in the proof
of Proposition \ref{ref-3.2.5-15} that $E=0$.

The case $p=0$ is similar.
\end{proof}
Now we employ our projective generators to construct a derived
equivalence. 
\begin{corollarys}
\label{ref-3.2.8-19} Assume that $\Pscr$ is a projective generator for
${}^p\Per(Y/X)$.  Put $A=\End_Y(\Pscr)$ and write ${}_A\Pscr$ to
emphasize the left $A$-structure on $\Pscr$
Then the functors
$ \RHom_Y({}_A\Pscr,-)$ and $-\Lotimes_A \,{}_A\!\Pscr$ define inverse equivalences between $D^b(\coh(Y))$ and
$D^b(A)$. These equivalences restrict to equivalences
between ${}^p\Per(Y/X)$ and $\mod(A)$.
\end{corollarys}
\begin{proof}
It is sufficient to show that the indicated functors 
define equivalences between ${}^p\Per(Y/X)$ and
$\mod(A)$ since the statement about derived categories then follows by
induction over triangles.  

Since we compute $ \RHom_Y({}_A\Pscr,-)$ with
injective resolutions in the second argument, we have that the
composition of $ \RHom_Y({}_A\Pscr,-)$ with the forgetful functor
$D^b(A)\r D^b(R)$ coincides with
$\RHom_Y({}_R\Pscr,-)$. So below we ignore the difference between
these functors.  

If $E\in {}^p\Per(Y/X)$ then by lemmas
\ref{ref-3.2.3-13} and \ref{ref-3.2.4-14} and Proposition
\ref{ref-3.2.6-17} we have $\Ext^i_Y(\Pscr,E)=0$ for $i\neq 0$. Hence on
${}^p\Per(Y/X)$ we have
$\RHom_Y(\Pscr,-)=\Hom_Y(\Pscr,-)$ and the latter sends
${}^p\Per(Y/X)$ to $\mod(R)$. 

Now we consider $-\Lotimes_A {}_A\Pscr$. Let $M\in
\mod(A)$. 
Breaking  a projective resolution of $M$ into short exact sequences and
using the exactness of $-\otimes_A \Pscr$ (in the notation of lemma
\ref{ref-3.2.1-11})  we easily deduce that for  $i\neq 0$:
${}^pH^i(M\Lotimes_A \Pscr)=0$  and ${}^pH^0(M\Lotimes_A
\Pscr)=M\otimes_A\Pscr$.
  Thus $-\Lotimes_A
{}_A\Pscr$ restricted to $\mod(A)$ coincides with $-\otimes_A
\Pscr$. Since $\Hom_Y(\Pscr,-)$ and $-\otimes_A
\Pscr$ define inverse equivalences between ${}^p\Per(Y/X)$ and
$\mod(A)$ we are done my lemma \ref{ref-3.2.1-11}
\end{proof}

As a final result in this section we show that in favourable cases
there is a relation between the projective objects in ${}^p\Per(Y/X)$
and $R$-Cohen-Macaulay modules.
\begin{lemmas}
\label{ref-3.2.9-20} Assume that $R$ is finitely generated over a field
  or else that it is a complete local ring containing a copy of its
  residue field. Besides our blanket hypotheses assume in addition
  that $X$ and $Y$ are Gorenstein of pure  dimension $n$ and
  that $f$ is birational and crepant (meaning:
  $f^\ast\omega_X=\omega_Y)$.  If $\Mscr$ is 
  a vector bundle on $Y$ then for every maximal ideal of $R$ we have
  $\depth_m \Gamma(Y,\Mscr)\ge n-1$.  If in addition 
$
H^1(Y,\Mscr)=H^1(Y,\Mscr^\ast)=0
$
then $\Gamma(Y,\Mscr)$ is a Cohen-Macaulay $R$-module.
\end{lemmas} 
\begin{proof} Let $\Dscr_X$, $\Dscr_Y$ be the Grothendieck dualizing
  complexes on $X$ and 
$Y$. One has $\Dscr_Y=f^!\Dscr_X$ and hence by the
Gorenstein/dimension hypotheses 
$\omega_Y=\Dscr_Y[-n]=f^!\Dscr_X[-n]=f^!\omega_X=f^!\Oscr_X\otimes_Y
f^\ast\omega_X$. Using the crepant hypothesis we obtain
$\Oscr_Y=f^!\Oscr_X$.  We compute for $\Mscr\in \coh(Y)$: $
\Ext^i_Y(\Mscr,\Oscr_Y)=\Ext^i_Y(\Mscr,f^!\Oscr_X)
=\Ext^i_R(R\Gamma(Y,\Mscr),R)
$. From this one easily deduces the assertions about
$\Gamma(Y,\Mscr)$. 
\end{proof}
\begin{propositions}
\label{ref-3.2.10-21}
Let the hypotheses be as in lemma \ref{ref-3.2.9-20} and assume
in addition that $X$ is normal. Then there
exists a Cohen-Macaulay module $M$ over $R$ such that if $A=\End_R(M)$
then $D^b(\coh(X))$ is derived equivalent to $D^b(A)$. In addition $A$
itself is Cohen-Macaulay. 
\end{propositions}
\begin{proof}
Let $\Mscr$ be a projective generator for
${}^{-1}\Per(Y/X)$. We have 
$H^1(Y,\Mscr^\ast)=0$ by definition and since $\Mscr$ is generated by
global sections we also have $H^1(Y,\Mscr)=0$. So $M=\Gamma(Y,\Mscr)$
is Cohen-Macaulay by lemma \ref{ref-3.2.9-20}. Let
$\Ascr=\uEnd_Y(\Mscr)$. As a vector bundle $\Ascr$ is self dual. We have
$\Ext^i_Y(\Mscr,\Mscr)=H^i(Y,\Ascr)$. Hence by lemmas
\ref{ref-3.2.3-13} and \ref{ref-3.2.9-20}
and Proposition \ref{ref-3.2.6-17} $A=\Gamma(Y,\Ascr)$ is
  Cohen-Macaulay, and hence reflexive since $\dim R\ge 2$.

There is an obvious map $A\r \End_R(M)$. This is an isomorphism
outside the exeptional locus in $X$ and this locus has codimension
$\ge 2$ in $X$. 
Since both objects are reflexive as
$R$-modules this implies that in fact $A=\End_R(M)$.
\end{proof}
Let $R$ be as in lemma \ref{ref-3.2.9-20}. If $A$ is an $R$-algebra
which is finite as an $R$-module then $A$ is said to be homologically
homogeneous if all simple $A$ modules have the same projective
dimension (necessarily $n$). It is shown in \cite{BH} that this
implies that $A$ has finite global dimension and is Cohen-Macaulay. In
fact the converse is also true and we leave it as a pleasant exercise
in homological algebra for the reader to check this.

The homological behaviour of homologically homogeneous rings closely
ressembles that of commutative rings of finite global dimension \cite{BH}.
\begin{corollarys} \label{ref-3.2.11-22} Let the hypotheses be as in
  Proposition 
  \ref{ref-3.2.10-21} but assume in addition that $Y$ is regular. Then
  there is a Cohen-Macaulay $R$-module $M$ such that
  $D^b(\coh(Y))$ is equivalent to $D^b(A)$ and such that $A=\End_R(M)$
  homologically 
  homogeneous,
\end{corollarys}
\begin{proof} We already  know that $A$ is Cohen-Macaulay. The fact
  that  $D^b(\coh(Y))$ is equivalent to $D^b(A)$ implies that $A$ has
  finite global dimension. Hence by the above discussion $A$ is
  homologically homogeneous.
\end{proof} 
 \subsection{General base}
\label{ref-3.3-23}
Now we assume that $X$ is a general noetherian scheme. We call an object
$\Pscr$ 
in ${}^p\Per(Y/X)$ a \emph{local projective generator} if 
there exists an open affine covering $X=\bigcup_i U_i$ of $X$ such
that for all $i$
 $\Pscr\mid f^{-1}(U_i)$ is a projective generator in
 ${}^p\Per(f^{-1}(U_i)/U_i)$.  
\begin{propositions}
\label{ref-3.3.1-24}
Assume that $\Pscr$ is a local projective generator for
${}^p\Per(Y/X)$. Put $\Ascr=f_\ast\uEnd_Y(\Pscr)$.  Then the functors
$Rf_\ast \uRHom_Y(\Pscr,-)$ and $f^{-1}(-)\Lotimes_{f^{-1}(\Ascr)}
\Pscr$ define inverse equivalences between $D^b(\coh(Y))$ and
$D^b(\coh(\Ascr))$. These equivalences restrict to equivalences
between ${}^p\Per(Y/X)$ and $\coh(\Ascr)$.
\end{propositions}
\begin{proof} This is easily proved by restricting to affine opens and
  invoking Corollary \ref{ref-3.2.8-19}.
\end{proof}
Under mild hypotheses we may construct a local projective generator
for ${}^p\Per(Y/X)$. 
\begin{propositions}
\label{ref-3.3.2-25}
Assume that $X$ is quasi-projective over a noetherian ring $S$. Then there
exists a local projective 
generator $\Pscr$ for 
${}^{-1}\Per(Y/X)$ such that $\Pscr^\ast$ is a local projective
generator for ${}^{-1}\Per(Y/X)$.
\end{propositions}
\begin{proof}
Let $\bar{X}$ be a projective $S$-scheme containing $X$ as an open
subset. Let $\bar{Y}$ be the closure of $Y$ under a locally closed
embedding $Y\r \PP^N_X\r \PP^N_{\bar{X}}$ and let $\bar{f}:\bar{Y}\r
\bar{X}$ be the corresponding projection morphism.

Let ${{\bar{\Lscr}}}$ be an $\bar{f}$-ample line bundle generated by global
sections 
 and choose an 
epimorphism $\phi:\bar{\Escr}\r R^1  f_\ast ({{\bar{\Lscr}}}^{-1})$ with
${{\bar{\Escr}}}=\Oscr_{\bar{X}}(-a)^b$.  By enlarging
$a$ we may assume 
\begin{equation}
\label{ref-3.2-26}
\Ext^i_{\bar{X}}({{\bar{\Escr}}},\bar{f}_\ast{{\bar{\Lscr}}}^{-1})=0\qquad
\text{for $i>0$}
\end{equation}
We have the following chain of equalities:
$
\Ext_{\bar{Y}}^1({{\bar{\Lscr}}},\bar{f}^\ast ({{\bar{\Escr}}}^\ast))
=\Ext_{\bar{Y}}^1(\bar{f}^\ast {{\bar{\Escr}}},{{\bar{\Lscr}}}^{-1})
=\Ext_{\bar{X}}^1({{\bar{\Escr}}},R\bar{f}_\ast{{\bar{\Lscr}}}^{-1})
=\Hom_{\bar{X}}({{\bar{\Escr}}},R^1\bar{f}_\ast{{\bar{\Lscr}}}^{-1})
$
where the last equality follows from \eqref{ref-3.2-26}. 

Hence the map $\phi$ provides us with an extension
\[
0\r  \bar{f}^\ast( {{\bar{\Escr}}})^\ast \r \bar{\Pscr}_0\r {{\bar{\Lscr}}}\r 0
\]
and restricting to $Y$ we obtain a corresponding extension
\[
0\r  f^\ast (\Escr)^\ast \r \Pscr_0\r \Lscr\r 0
\]
and it is clear from the construction that on small open affines the
extension restricts to \eqref{ref-3.1-16}. Thus $\Pscr=\Pscr_0\oplus \Oscr_Y$
is our required local projective generator.
\end{proof}

\subsection{The formal case}
\label{ref-3.4-27}
We now discuss the case $X=\Spec R$ with $R$ a noetherian complete
local ring with maximal ideal $m$ such that $k=R/m$ is algebraically
closed and $k\subset R$ (at the cost of some technicalities one can
get by with less). It is not necessary to assume $\operatorname{char} k=0$. The formal
case is interesting since it contains some features not present in the
general case. Also our main result in \S\ref{ref-4.4-58} is proved by
reduction to the formal case.

Let $x$ be 
the unique closed point of $X$ (i.e the defining ideal of $x$ is $m$)
and let $C=f^{-1}(x)$. $C$ is either scheme-theoretically a point (if
$f$ is an isomorphism) or else it has dimension one.

First we discuss the structure of $C_{\red}$. This is of course
well-known.
\begin{lemmas} \label{ref-3.4.1-28} $C_{\red}$ is either a point or a
  tree of $\PP^1$'s 
  with normal crossings.
\end{lemmas}
\begin{proof} 
Assume that $C$ is not a point.
Since $C$ is connected, it is clear that the only global sections of
$\Oscr_{C_{\red}}$ are scalar. Furthermore since $\Oscr_{C_{\red}}$ is a quotient of
$\Oscr_Y$ we deduce that
$H^1(C_{\red},\Oscr_{C_{\red}})=0$. Thus $C_{\text{red}}$ is a reduced
projective curve of arithmetic genus zero. It is well-known that this must be a
tree of $\PP^1$'s.
\end{proof}
We also have the following.
\begin{lemmas}
\label{ref-3.4.2-29}
$C$ is Cohen-Macaulay (i.e. $C$ has no embedded components) and
$H^0(C,\Oscr_C)=k$. 
\end{lemmas}
\begin{proof} This clear if $C$ is a point. Assume that this is not
  the case.
We have an exact sequence
\[
0\r m\Oscr_Y \r \Oscr_Y \r \Oscr_C\r 0
\]
Since $m\Oscr_Y$ is generated by global sections we have
$H^1(Y,m\Oscr_Y)=0$. Hence $H^0(C,\Oscr_C)$ is a quotient of
$H^0(Y,\Oscr_Y)$ by an ideal containing $m$. It follows
that $H^0(C,\Oscr_C)=k$. Any embedded component of $C$ would be
zero-dimensional and hence would give rise to extra sections. So such
embedded components cannot exist. 
\end{proof}
We will denote the irreducible components of $C$
by $(C_i)_{i=1,\ldots,n}$. If $C$ is a point then we set $n=0$. 

In each of the $C_i$ we choose a point
$y_i$ not lying on the other components and not lying on the possible
finite number of embedded zero dimensional components of $Y$.

The following result is well-known.
\begin{lemmas} \label{ref-3.4.3-30} The map $\Lscr\mapsto \deg(\Lscr\mid C_i)_{i=1,\ldots,n}$ defines
  an isomorphism
$\Pic(Y)\cong \ZZ^n$.
\end{lemmas}
\begin{proof}
Let $\hat{Y}$ be the formal scheme associated to $(Y,C)$.  We have
$H^1(\hat{Y},\Oscr_{\hat{Y}})=H^1(Y,\Oscr_Y)=0$ \cite[\S 4]{EGAIII}
and by Grothendieck's existence theorem  we also have
$\Pic(X)=\Pic(\hat{X})$. 
Consider the exact sequence
\begin{equation}
\label{ref-3.3-31}
0\r m\Oscr_{\hat{Y}}\r \Oscr^\ast_{\hat{Y}} \r \Oscr^\ast_C \r 0
\end{equation}
As in the proof of the previous lemma we have
$H^1(\hat{Y},m\Oscr_{\hat{Y}})=0$. Hence from the long exact sequence associated to
\eqref{ref-3.3-31} we deduce $\Pic(\hat{Y})=\Pic(C)$. Now $C$ is
either a point or a
projective curve with  $H^1(C,\Oscr_C)=0$. It is well known that this
implies $\Pic(C)=\ZZ^n$. It is easy to see that the isomorphism has
the indicated form. 
\end{proof}
The previous lemma implies the existence of line bundles $\Lscr_i$ on
$Y$ such that $\deg(\Lscr_i\mid C_j)=\delta_{ij}$.
One may improve slightly on this fact.
\begin{lemmas} \label{ref-3.4.4-32} There exist  a closed subscheme $D_i$ in $Y$  with
  $\Oscr_Y(-D_i)$ invertible   such that scheme theoretically we have
\[
D_i\cap C_j=
\begin{cases}
\{y_i\}&\text{if $i=j$}\\
\emptyset &\text{otherwise}
\end{cases}
\]
\end{lemmas}
\begin{proof} If $D$ is an arbitrary closed subscheme of $Y$ then 
the connected components of $D$ coincide
  with the connected components of $D\cap C$. 

It is easy to see that on an affine neighborhood $U_i$ of $y_i$ we can
choose a non-zero divisor $z\in \Gamma(U_i,\Oscr_{U_i})$ such that
$V(z)\cap C_i=\{y_i\}$ scheme-theoretically. Let $D'$ be the closure
of $V(z)$ in $Y$. We let $D_i$ be the component of $D'$ containing $y_i$.
\end{proof}
Let $\Pic^+(Y)$ and $\Pic^{++}(Y)$ be the isomorphism classes
of line bundles which are respectively
generated by gobal sections and ample. The following
is well-known:
\begin{lemmas} We have
\begin{align} \label{ref-3.4-33}
\Pic^+(Y)&=\{ \Lscr\in \Pic(Y)\mid \deg(\Lscr\mid C_i)\ge 0\text{ for
  all $i$}\}\\
\Pic^{++}(Y)&=\{ \Lscr\in \Pic(Y)\mid \deg(\Lscr\mid C_i)> 0\text{ for
  all $i$}\}
\end{align}
\end{lemmas}
\begin{proof} The statement about $\Pic^{++}(Y)$ follows from
  \cite[Prop. 2.7]{Ke1}. For the statement about $\Pic^+(Y)$ we use
  the fact that if $\Lscr$ is as in \eqref{ref-3.4-33} then by lemma \ref{ref-3.4.4-32}
  $\Lscr=\Oscr_Y(E)$ with $E=\sum n_i D_i$,  $n_i\ge 0$. The support
  of $\Oscr_E(E)$ is finite over $X$ and hence affine. From the exact
  sequence
\[
0\r \Oscr_Y \r \Oscr_Y(E)\r \Oscr_E (E)\r 0
\]
and the fact that $H^1(Y,\Oscr_Y)=0$ we obtain that $\Oscr_Y(E)$ is
generated by global sections.
\end{proof}
\subsection{Artin-Verdier theory}
\label{ref-3.5-34}
In \cite{AV} Artin and Verdier show that non-trivial indecomposable
Cohen-Macaulay 
modules on a rational double point are in one-one correspondence with
the irreducible components of the exceptional divisor in a
minimal resolution. In this section we show that part of that theory
generalizes to the setting introduced in \S\ref{ref-3.4-27}. 
\begin{lemmas}
\label{ref-3.5.1-35}
Let $\Mscr$ be a vector bundle of rank $r$ on $Y$ generated by global
sections. Then $\Mscr$ occurs in exact sequences of the form
\begin{gather}
\label{ref-3.6-36}
0\r \Oscr^{r-1}_Y \r \Mscr \r \Lscr\r 0\\
\label{ref-3.7-37}
0\r \Lscr^{-1}\r \Oscr_Y^{r+1} \r  \Mscr\r 0 
\end{gather}
where $\Lscr=c_1(\Mscr)$.
\end{lemmas}
\begin{proof}
This is proved in a similar way as  \cite[Lemma 1.2,1.5]{AV}. The case
where $C$ is a point is clear so we assume this is not the case.
 
Let $\phi: \Oscr^{r-1}_Y \r \Mscr$, $\theta:\Oscr_Y^{r+1} \r  \Mscr$
be obtained from choosing respectively $r-1$ and $r+1$ generic
sections of $\Mscr$.

Let $J_i$ be the ideal sheaf of $C_i$. From the exact sequence
\[
0\r J_i \r \Oscr_Y\r \Oscr_{C_i}\r 0
\]
we deduce that $H^1(Y,J_i)=0$. Since $\Mscr$ is generated by global
sections we also have $H^1(Y,J_i\Mscr )=0$. It follows that
$H^0(Y,\Mscr)\r H^0(Y,\Mscr\otimes_Y \Oscr_{C_i})$ is surjective for
all $i$. In particular generic sections of $\Mscr$ correspond to
generic sections of $H^0(Y,\Mscr\otimes_Y \Oscr_{C_i})$.

We claim that $\phi$ and $\theta$ have maximal rank everywhere. This
may be checked on closed points and hence it may be checked by
restricting to $C_i$. By lemma \ref{ref-3.4.1-28} $C_i=\PP^1$ and on $\PP^1$
the statement is elementary. 

Hence $\Lscr=\coker \phi$ and $\Lscr'=\ker \theta$ are line
bundles. That $\Lscr=\Lscr^{\prime -1}=c_1(\Mscr)$ is clear.
\end{proof}

\begin{lemmas}
\label{ref-3.5.2-38}
The vector bundles $\Mscr$ in $\frak{V}$ are precisely the ones which occur
as middle term in an exact
sequence  with $\Lscr\in \Pic^+(Y)$
\begin{equation}
\label{ref-3.8-39}
0\r \Oscr_Y^{r-1} \r \Mscr \r \Lscr \r 0
\end{equation}
which is determined by a set of $r-1$ generators  of
$H^1(Y,\Lscr^{-1})$. $\Mscr$ is uniquely determined by
$\Lscr=c_1(\Mscr)$, up to 
the addition of copies of $\Oscr_Y$. 

The vector bundles $\Nscr$ in $\frak{V}^\ast$ are precisely the ones
which occur in an exact sequence with $\Lscr\in \Pic^+(Y)$
\begin{equation}
\label{ref-3.9-40}
0\r \Nscr \r \Oscr_Y^{r+1} \r \Lscr \r 0
\end{equation}
which is determined by a set of $r+1$ generators  of
$H^0(Y,\Lscr)$. $\Nscr$ is uniquely determined by
$\Lscr=c_1(\Nscr)^{-1}$, up to 
the addition of copies of $\Oscr_Y$. 
\end{lemmas}
\begin{proof} 
That every $\Mscr$ which occurs in an exact sequence
\eqref{ref-3.8-39} is in $\frak{V}$ is clear so we concentrate on
the converse. According to lemma \ref{ref-3.5.1-35} $\Mscr$ occurs in
  an exact sequence 
\[
0\r \Oscr_Y^{r-1} \r \Mscr \r \Lscr\r 0
\]
and the fact that $H^1(Y,\Mscr^\ast)=0$ implies that this exact
sequence is determined by a set of $r-1$-generators of $H^1(Y,\Lscr^{-1})$.

Any set of generators $H^1(Y,\Lscr^{-1})$ contains a minimal set of
generators. Adding extra generators corresponds to adding free
summands to $\Mscr$. It is easy to see that two minimal sets of
generators yield isomorphic $\Mscr$. 

The proof for $\Nscr$ is similar (\eqref{ref-3.9-40} is obtained
by dualizing \eqref{ref-3.7-37}). 
\end{proof}
 From the previous lemmas we extract the following corollary.
\begin{propositions} 
\label{ref-3.5.3-41} The map
\[
\phi:\frak{V}\r \ZZ\times \Pic^+(Y):\Mscr\r (\rk(\Mscr),c_1(\Mscr))
\]
 is an
  injection on isomorphism classes of objects.
\end{propositions}
In order to understand the image of $\phi$ we have to investigate the
indecomposable objects in $\frak{V}$.

Let $\Mscr_i\in \frak{V}$ be the extension
\begin{equation}
\label{ref-3.10-42}
0\r \Oscr_Y^{r_i-1}\r \Mscr_i \r \Lscr_i\r 0
\end{equation}
associated to a minimal set $r_i-1$ of generators of
$H^1(Y,\Lscr^{-1}_i)$ where $\Lscr_i$ is as in
\S\ref{ref-3.4-27}. 
 We also define $\Mscr_0=\Oscr_Y$. 
\begin{propositions}
\label{ref-3.5.4-43}
The $\Mscr_i$ are indecomposable and furthermore every indecomposable
object in $\frak{V}$ is isomorphic to one of the 
$\Mscr_i$. If $i>0$ then $\rk \Mscr_i$ is equal to the multiplicity of $C_i$
in $C$. 
\end{propositions}
\begin{proof} In the proof below we assume $i>0$. 
It is clear from Proposition \ref{ref-3.5.3-41} that if $\Mscr_i$
  is decomposable then it is of the form $\Oscr_Y^a\oplus \Mscr'$ with
  $\Mscr'$ indecomposable (this corresponds to the only possible way of
  decomposing $(\rk(\Mscr_i),c_i(\Mscr_i))$ in $\ZZ\times \Pic^+(Y)$).
 Let $r'$ be the rank
  of $\Mscr'$. Thus $r'\le r_i$. According to lemmas \ref{ref-3.5.2-38}
  $\Mscr'$ appears in an exact sequence
\[
0\r \Oscr^{r'-1}_Y\r \Mscr'\r \Lscr_i\r 0
\]
corresponding to a set of $r'-1$ generators of
$H^1(Y,\Lscr_i^{-1})$. But the minimal number of generators of
$H^1(Y,\Lscr_i^{-1})$ is $r_i-1$. Thus $r'\ge r_i$ and hence
$r'=r_i$. Thus $\Mscr_i=\Mscr_i'$. In other words: $\Mscr_i$ is
indecomposable. 

Now let $\Mscr$ be an arbitrary object in $\frak{V}$. We may find an
object $\Rscr'=\Mscr_1^{\oplus a_1}\oplus \cdots \oplus \Mscr_n^{\oplus
  a_n}$ such that $c_1(\Mscr)=c_1(\Rscr')$. So by Proposition
\ref{ref-3.5.3-41} $\Mscr$ and $\Rscr'$ differ only by copies of
$\Oscr_Y$. Since $\Rscr'$ can have no summands isomorphic to $\Oscr_Y$
these must be in $\Mscr$. So there is some $\Rscr=\Oscr_Y^{\oplus
  a_0}\oplus \Rscr'$ such that $\Mscr\cong \Rscr$. This shows that
every indecomposable object in $\frak{V}$ is isomorphic to one of the
$\Mscr_i$.

To prove the statement on the rank of $\Mscr_i$ we need to compute the
minimal number of generators of  $H^1(Y,\Lscr_i^{-1})$. By lemma
\ref{ref-3.4.3-30} we have $\Lscr_i\cong\Oscr_Y(-D_i)$. From the
exact sequence
\[
0\r \Oscr_Y(-D_i)\r \Oscr_Y\r \Oscr_{D_i}\r 0
\]
we obtain an exact sequence
\[
R\r \Gamma(Y,\Oscr_{D_i})\r H^1(Y,\Oscr_Y(-D_i))\r 0
\]
The first map is a ring map so its image cannot be in
$m\Gamma(Y,\Oscr_{D_i})$. Thus the minimal number of generators of
$H^1(Y,\Oscr_Y(-D_i))$ is one less than the minimal number of
generators of $\Gamma(Y,\Oscr_{D_i})$. The latter is equal to the
dimension of the vector space $\Gamma(Y,\Oscr_{D_i}/m\Oscr_{D_i})$
which is precisely $C\cdot D_i$. Given the definition of $D_i$ this is
also the multiplicity of $C_i$ in $C$.
\end{proof}
In the sequel we put $\Nscr_i=\Mscr_i^\ast$. By lemma \ref{ref-3.5.2-38}
$\Nscr_i$ occurs in an exact sequence 
\begin{equation}
\label{ref-3.11-44}
0\r \Nscr_i \r \Oscr^{r_i+1}_Y \r \Lscr_i \r 0
\end{equation}
\begin{theorems}
\label{ref-3.5.5-45}
The indecomposable projective objects in ${}^{-1}\Per(Y/X)$ are
the $(\Mscr_i)_i$. The indecomposable projective objects in
${}^0\Per(Y/X)$ are the $(\Nscr_i)_i$. The projective generators in
${}^{-1}\Per(Y/X)$ are of the form $\oplus_i \Mscr_i^{\oplus a_i}$,
$a_i>0$. The projective generators in
${}^0\Per(Y/X)$ are of the form $\oplus_i \Nscr_i^{\oplus b_i}$,
$b_i>0$. 
\end{theorems}
\begin{proof} This is a combination of Propositions
  \ref{ref-3.2.6-17}, \ref{ref-3.2.7-18} and \ref{ref-3.5.4-43}.
\end{proof}
Summarizing  we have proved the following result. 
\begin{theorems} 
\label{ref-3.5.6-46} We have ${}^{-1}\Per(Y/X)\cong \mod(A)$,
${}^0\Per(Y/X)\cong \mod(A^\circ)$ where $A$ is  finitely generated
as a module over $R$ and $A/\rad(A)\cong k^{n+1}$.
\end{theorems}
\begin{proof} It suffices to take as projective generator $\oplus_i
  \Mscr_i$ for ${}^{-1}\Per(Y/X)$ and  $\oplus_i
  \Nscr_i$ for ${}^0\Per(Y/X)$
\end{proof}
It is possible to give explicitly the $n+1$ simple objects in
${}^p\Per(Y/X)$. This was independently observed by Bridgeland. This
result will not be used afterwards.
\begin{propositions} The $n+1$ simple objects in ${}^{-1}\Per(Y/X)$
  are $\Sscr_0=\Oscr_C$ and $\Sscr_i=\Oscr_{C_i}(-1)[1]$, $i=1,\ldots,n$.
\end{propositions}
\begin{proof}  By Theorem \ref{ref-3.5.6-46} it
  is sufficient to prove 
  $\Hom(\Mscr_i,\Sscr_j)=\delta_{ij}\cdot k$. 

The case $i=0$ is clear so we assume $i\neq 0$. If $j\neq 0$ then we
have $\Hom_Y(\Mscr_i,\Oscr_{C_j})=\Ext^1_Y(\Lscr_i,\Oscr_{C_j})=
\Ext^1_{C_j}(\Lscr_i\mid C_j,\Oscr_{C_j})$. Now $\Lscr_i\mid C_j$ is
either $\Oscr_{C_j}(1)$ or $\Oscr_{C_j}$ depending on whether $i=j$ or
$i\neq j$. This finishes the proof in the case $i,j \neq 0$.

The only
remaining case is $i\neq 0$, $j=0$. We need to compute
$\Hom_Y(\Mscr_i,\Oscr_C)=\Hom_C(\Mscr_i\otimes_Y \Oscr_C,\Oscr_C)$ for
$i\neq 0$.  From the right exactness of $H^1(Y,-)$ we deduce  for
$M\in \mod(R)$: $H^1(Y,\Lscr_i^{-1}\otimes_R
M)=H^1(Y,\Lscr_i^{-1})\otimes_R M$. Applying this with $M=R/m$ we see
that a set of minimal generators of $H^1(Y,\Lscr_i^{-1})$ as $R$-module yield
a basis of  $H^1(Y,\Lscr_i^{-1}\otimes_Y \Oscr_C)$ as $k$-vector space.
This implies that in the
long exact 
sequence for $\Hom_Y(-,\Oscr_C)$ applied to \eqref{ref-3.10-42}, the
first connecting map is an isomorphism.  Thus we obtain
$\Hom_Y(\Mscr, \Oscr_C)=\Gamma(C,\Lscr_i^{-1}\otimes_Y \Oscr_C)$. Now
$\Gamma(C,\Lscr_i^{-1}\otimes_Y \Oscr_C)=\Gamma(C,\Oscr_C(-D_i))$ and
this is zero 
since $\Gamma(C,\Oscr_C(-D_i)\subset \Gamma(C,\Oscr_C)=k$, and the
generator $1$ for  $\Gamma(C,\Oscr_C)$ is nowhere vanishing.
\end{proof}
The following result is proved in a similar way.
\begin{propositions} The simple objects in ${}^0\Per(Y/X)$ are
  $\omega_C[1]$ and $\Oscr_{C_i}(-1)$. If $(-)^D$ denotes the
  Grothendieck duality functor on $D^b(\coh(Y))$ \cite{RD} then
  $(-)^D$ defines a duality between the categories of finite length
  objects in ${}^p\Per(Y/X)$ and ${}^{-1-p}\Per(Y/X)$. 
\end{propositions}

\section{Flops}
\subsection{Some special Cohen-Macaulay modules}
\label{ref-4.1-47}
Below $P$ is a $n\ge 3$-dimensional regular complete local ring and
$R$ is a normal integral Gorenstein ring which is of rank two over $P$. In
particular $R$ is free over $P$. By $\Cl(R)$ we denote the \emph{class
  group} of $R$ \cite{Fossum}. Its elements are isomorphism classes
of rank one reflexive $R$-modules and the multiplication is given by $I\cdot
J=(I\otimes_R J)^{\ast\ast}$.

Let $I$ be a reflexive rank one module over $P$ of depth $\ge
n-1$.  Thus $\Ext^i_R(I,R)=0$ for
$i\ge 2$. 

We consider exact sequences of $R$-modules of the
form 
\begin{gather}
\label{ref-4.1-48}
0\r R^{r-1}\r M\r I \r 0\\
\label{ref-4.2-49}
0\r N\r R^{s+1}\r I \r 0
\end{gather}
where we assume that \eqref{ref-4.1-48} is given by a set of $r-1$ generators
of $\Ext^1_R(I,R)$. 
Looking at the long exact sequences for
$\Hom_R(-,R)$ we find $\Ext^i_R(M,R)=\Ext^i_R(N,R)=0$ for $i>0$. Hence
$M,N$ are (maximal) Cohen-Macaulay modules. 

Now $M,N$ are determined by $I$ up to adding of free summands. We
denote by $M(I)$, $N(I)$  the Cohen-Macaulay modules obtained
from $M$, $N$ by deleting the free summands. 
\begin{remarks}  Indecomposable Cohen-Macaulay modules are not always
  of the form 
  $M(I)$, $N(I)$.  Recall that 
  if $R$ is Cohen-Macaulay but non regular there exist always
  non-free indecomposable Cohen-Macaulay modules (for example given by
  a summand of a suitable 
  syzygy of $R/m$). But if  $R$ is in addition factorial
  then the only 
  Cohen-Macaulay module of the form $M(I)$, $N(I)$  are trivial which
  yields a contradiction. A concrete counter example is given by
  $k[[x,y,z,t]]/(x^2+y^2+z^2+t^3)$.

  The author thinks it is an intriguing question to understand
  precisely which indecomposable Cohen-Macaulay modules are of the
  form $M(I)$ or $N(I)$. It appears possible to solve this problem for
  complete local rings with a three-dimensional isolated Gorenstein
  terminal singularity. We will come back on this in a subsequent
  paper. 
\end{remarks}

In the sequel we need the following technical result.
\begin{propositions} \label{ref-4.1.2-50} We have $N(I)\cong M(I^{-1})$.
\end{propositions} 
\begin{proof}  Note that for an $R$-module $M$ we have $\depth_R
  M=\depth_P M$ and furthermore $M$ is reflexive as $R$ module if and
  only if it is reflexive as $P$-module. 
 
Hence since $R$
  is free of rank two over $P$, $R\otimes_P I$ is reflexive and has
  depth $\ge n-1$. Its rank over $R$ is two. 

Since $P$ is regular local the projective dimension of $I$ over $P$ is
one. Since $R$ is flat as $P$ module it follows that $R\otimes_P I$
also has projective dimension one.

Let $K$ be the kernel of $R\otimes_P I\r I$. 
$K$ is also reflexive of rank one and it has depth $\ge n-1$. To
compute  $K$  we note 
that if we assign to a reflexive $R$-module $M$ the element of
$\Cl(R)$ 
given by $(\wedge^{\rk M} M)^{\ast\ast}$  then
$c_1$ is multiplicative on short exact sequences.  

Applying this to the exact sequence
\[
0\r K \r R\otimes_P I \r I \r 0
\]
and using the fact that $c(R\otimes_P I)=R$ since $R\otimes_P I$ has
finite projective dimension, we find $K\cong I^{-1}$. 

Now suppose that we have an exact sequence as in \eqref{ref-4.2-49}. After
possibly adding free summands to $N$ and $R^{s+1}$ we may construct a
diagram
\[
\begin{CD}
0 @>>> N @>>> R^{s+1} @>>> I @>>> 0 \\
@. @. @VVV @|\\
@. @. R\otimes_P I @>>> I @>>> 0\\
@. @. @VVV\\
@. @. 0
\end{CD}
\]
Using the fact that $R\otimes_P I$ has projective dimension one 
this diagram may be completed to
\[
\begin{CD}
@. 0 @. 0\\
@. @VVV @VVV\\
@. R^{s-1} @= R^{s-1}\\
@. @VVV @VVV\\
0 @>>> N @>>> R^{s+1} @>>> I @>>> 0 \\
@. @VVV @VVV @|\\
0 @>>> I^{-1} @>>> R\otimes_P I @>>> I @>>> 0\\
@. @VVV @VVV\\
@. 0 @. 0
\end{CD}
\]
By construction $N$ is Cohen-Macaulay so $\Ext^1_R(N,R)=0$. It follows
that the left most vertical exact sequence is obtained from a set of
generators of
$\Ext^1_R(I^{-1},R)$. So up to free summands we have $N(I)=N=M(I^{-1})$.
\end{proof}

\subsection{Application to hypersurface singularities.}
\label{ref-4.2-51}

The following lemma will be used.
\begin{lemmas}
Let $f:Y\r X$ be a projective birational map between normal 
noetherian schemes such
that the exceptional 
locus of $f$ has codimension $\ge 2$ in $Y$. 
\label{ref-4.2.1-52} The functor $f_\ast$ restricts to an equivalence
  between the category of reflexive $\Oscr_Y$-modules and the category of
  reflexive $R$-modules.
\end{lemmas}
\begin{proof} We may assume that $X=\Spec R$ is affine. Let $\Mscr$ be a reflexive $\Oscr_Y$ module. We need to
  show that for every prime ideal $p$ in the exceptional locus of $f$ in
  $X$ we have 
  $\Hom_X(p,f_\ast\Mscr)=\Gamma(Y,\Mscr)$. Now
  $\Hom_X(p,f_\ast\Mscr)=\Hom_Y(f^\ast
  p,\Mscr)$. By hypotheses the exceptional locus of $f$ has codimension $2$,
  and hence $f^\ast p$ is equal to $\Oscr_Y$ in codimension one. Since  $\Mscr$
  is reflexive this implies  $\Hom_Y(f^\ast
  p,\Mscr)=\Hom_Y(\Oscr_Y,\Mscr)=\Gamma(Y,\Mscr)$.

That $f_\ast$ yields an equivalence between the category of
reflexive $\Oscr_Y$-modules and the category of reflexive $R$-modules
also follows from the hypotheses on the codimension of the
exceptional locus.
\end{proof}

Let $\Cl(Y)$ be the group of Weil divisors on $Y$ modulo the principal
divisors. The previous lemma gives us a canonical identification
between $\Cl(Y)$ and $\Cl(R)$. 

We  specialize the situation of  \S\ref{ref-3.4-27}. Now 
 $R$  is a
normal complete local $k$-algebra of dimension $n \ge 3$ with residue
 field $k$ and with a
 canonical hypersurface
 singularity of 
 multiplicity two. According to \cite[Cor 5.24]{KM} $X$ has rational
 singularities. 

Let $f:Y\r X$
be 
a birational projective map such that $Y$ is normal 
Gorenstein. Assume in addition that the exceptional locus of $f$ has
codimension $\ge 2$ in $Y$. It is easy to see that $Y$ also has canonical and
hence rational singularities (again by \cite[Cor 5.24]{KM})
and therefore $Rf_\ast \Oscr_Y=\Oscr_X$

Note that by lemma \ref{ref-3.2.9-20}
$\Pic(Y)$ under the identification $\Cl(Y)\cong \Cl(R)$ is sent to
reflexive ideals of 
depth $\ge n-1$.

Let $\Mscr_i$, $\Nscr_i$ be as in the previous section. Then we deduce
 from lemma \ref{ref-3.2.9-20} that
$\Gamma(Y,\Mscr_i)$ and $\Gamma(Y,\Nscr_i)$ are Cohen-Macaulay.

We let $\Lscr_i$ be as in \S\ref{ref-3.4-27} and we put
$I_i=\Gamma(Y,\Lscr)$. Thus the $I_i$ are reflexive $R$-modules of
rank one and depth $\ge n-1$.
Put $M_i=M(I_i)$, $N_i=N(I_i)$ for $i>0$ and $M_0=N_0=R$.

The following observation is crucial.
\begin{lemmas}
\label{ref-4.2.2-53}
We have $\Gamma(Y,\Mscr_i)=M_i$, $\Gamma(Y,\Nscr_i)=N_i$. 
\end{lemmas}
\begin{proof} The case $i=0$ is trivial so we assume $i>0$.
We consider the first equality. Applying $\Gamma(Y,-)$ to \eqref{ref-3.10-42}
we obtain an 
exact sequence
\begin{equation}
\label{ref-4.3-54}
0\r R^{r_i-1} \r  \Gamma(Y,\Mscr_i) \r I_i \r 0
\end{equation}
According to Proposition \ref{ref-3.5.4-43}, lemma \ref{ref-4.2.1-52}
 and the above discussion we have  that $\Gamma(Y,\Mscr_i)$ 
is an  indecomposable Cohen-Macaulay $R$-module. Applying $\Hom(-,R)$ to
\eqref{ref-4.3-54} we see that \eqref{ref-4.3-54} is associated to a set
$r_i-1$ generators of $\Ext^1_R(I_i,R)$. Thus $\Gamma(Y,\Mscr_i)$
is obtained from $M_i$ by adding free summands. Since
$\Gamma(Y,\Mscr_i)$ is indecomposable we obtain
$\Gamma(Y,\Mscr_i)=M_i$. 

Now we consider the second equality. Applying $\Gamma(Y,-)$ to \eqref{ref-3.11-44}
 we obtain an 
exact sequence
\begin{equation}
\label{ref-4.4-55}
0\r \Gamma(Y,\Nscr_i)\r R^{r_i+1}  \r I_i \r 0
\end{equation}
Now since $\Nscr_i=\Mscr_i^\ast$
 it follows from  Proposition \ref{ref-3.5.4-43}, lemma \ref{ref-4.2.1-52},
 and the above discussion that $\Gamma(Y,\Nscr_i)$ 
is an indecomposable Cohen-Macaulay $R$-module. 

Thus $\Gamma(Y,\Nscr_i)$
is obtained from $N_i$ by adding free summands. Since
$\Gamma(Y,\Nscr_i)$ is indecomposable we obtain
$\Gamma(Y,\Nscr_i)=N_i$. 
\end{proof}
\begin{remarks} It follows from the above proof that the $M_i$, $N_i$
  are indecomposable and furthermore that they occur in exact
  sequences 
\begin{gather*}
0\r R^{r_i-1} \r M_i \r I_i \r 0\\
0\r N_i \r R^{r_i+1} \r I_i \r 0
\end{gather*}
In addition lemma \ref{ref-4.2.1-52} together with the fact that
$\Nscr_i=\Mscr_i^\ast$ implies  $N_i=M^\ast_i$.
\end{remarks}
\subsection{Formal flops}
\label{ref-4.3-56}
We recycle the notations of the previous section. According to
\cite[Prop. 2.3]{Kollar} the map $f:Y\r X$ has a flop $f^+:Y^+\r X$. 

The construction of $f^+$ is as
follows. Write $R=k[[x_1,x_2,\ldots,x_{n+1}]]/(x_1^2+f(x_2,\ldots,x_{n+1}))$
and let $\sigma:X\r X$ be given by $(x_1,x_2,\ldots,x_{n+1})\r
(-x_1,x_2,\ldots,x_{n+1})$. Then $Y^+=Y$ and $f^+=\sigma\circ f$. 

Below we use the same notations for $Y$ as for $Y^+$ but we adorn the
latter with a superscript ``$+$''.  To fix the numbering of the
$C_i^+$ 
we use the identification  $Y^+=Y$
to put $C^+_i=C_i$, $D^+_i=D_i$. 

 Lemma \ref{ref-4.2.1-52} gives us a canonical
identification $\Cl(Y)=\Cl(R)=\Cl(Y^+)$. Since $\sigma$ induces the
operation $I\mapsto I^{-1}$ on $\Cl(R)$ (see \cite[Example 2.3]{Kollar}) we
find that under this identification $\Lscr^+_i$ on $Y^+$ is equivalent to
$\Lscr^{-1}_i$ on $Y$. In particular $I^+_i=I^{-1}_i$. 

We now 
obtain our key result.
\begin{propositions} 
\label{ref-4.3.1-57} One has
$M_i^+\cong N_i$, $N_i^+\cong M_i$. 
\end{propositions}
\begin{proof}
This is a direct application of Proposition \ref{ref-4.1.2-50}. Indeed we
have $M^+_i=M(I_i^+)=M(I_i^{-1})=N(I_i)=N_i$. The argumentation for
$N^+_i$ is identical.
\end{proof}
\subsection{Global flops}
\label{ref-4.4-58}
\begin{lemmas} 
Let $f:Y\r X$ be a projective birational map between normal 
varieties over $k$ of dimension $n\ge 3$ such
that the exceptional 
locus of $f$ has codimension $\ge 2$ in $Y$. Assume that $X$ has 
hypersurface singularities (not necessarily isolated) of multiplicity $\le 2$. 

Under these conditions the flop of $f$ exists and is unique. More
precisely there exists a unique morphism 
$f^+:Y^+\r X$ with the following properties
\begin{enumerate}
\item $f^+$ is projective and birational and $Y^+$ is a normal 
  variety. The maximal dimensions of the fibers of $f$ and $f^+$ are
  the same.
  The exceptional locus of $f^+$ has codimension $2$ in $Y^+$. If $Y$
  is Gorenstein then so is $Y^+$.
\item Under the identifications obtained from lemma \ref{ref-4.2.1-52}:
$
\Cl(Y)\cong \Cl(X) \cong \Cl(Y^+)
$, 
$\Pic(Y)$ corresponds to $\Pic(Y^+)$. 
\item If $E$ is an $f$-nef (resp.
$f$-ample) divisor on $Y$ then $-E$ corresponds to an $f$-nef (resp. $f$-ample)
divisor on  $Y^+$.
\end{enumerate}
\end{lemmas}
\begin{proof}  
Choose an $f$-ample  Cartier divisor $D$ on $Y$ and
  identify it with a Weil divisor on $X$, also denoted by
  $D$.  Then $Y^+$, if it exists, is the $D$-flop of $f$ (in the
  terminology of \cite[\S 2]{Kollar}). Hence $Y^+$ is unique if it exists.

According to \cite[Cor. 6.4]{KM} for $Y^+$ to exist at least
  the sheaf of graded ring $\Sscr=\oplus_n \Oscr_X(-nD)$ should be a
  sheaf of finitely generated $\Oscr_X$-algebras.  In that case
  $Y^+=\underline{\Proj} \,\Sscr$.

According to
  \cite[Prop. 6.6, Cor. 6.7]{KM} the finite generation of $\Sscr$  can
  be checked in the completions of the closed points of $X$.

In the formal case we have $X=\Spec R$ where $R$ is a hypersurface
singularity of multiplicity one or two.
In the case of multiplicity one there is nothing to prove  and in the
case 
of multiplicity two the flop exists by
\cite[Example 2.3]{Kollar}. From the construction of $f^+$ (see
\S\ref{ref-4.3-56}) 
it easily follows that $f^+:Y^+\r X$ has the properties
(1)(2)(3) listed in the statement of the theorem.  Since these
properties may also be checked in the completions of the closed points
they hold in the general  case.
\end{proof}
Now we give our main theorem.
\begin{theorems}
\label{ref-4.4.2-59}
Let $f:Y\r X$ be a projective birational map between normal
quasi-projective  Gorenstein
$k$-varieties of dimension $n\ge 3$ with fibers of dimension $\le 1$ and
assume that the exceptional 
locus of $f$ has codimension $\ge 2$. Assume that $X$ has canonical
hypersurface singularities of multiplicity $\le 2$.  Let $f^+:Y^+\r X$
be the flop of $f$.  Then
  $D^b(\coh(Y))$ and
  $D^b(\coh(Y^+))$ are equivalent and we may choose  this equivalence 
in such a way that ${}^{-1}\Per(Y/X)$ corresponds to
 ${}^{0}\Per(Y^+/X)$.
\end{theorems}
\begin{proof}
 Let $\Pscr$ be a local projective generator of
 ${}^{-1}\Per(Y/X)$ ($\Pscr$ exists by Proposition
 \ref{ref-3.3.2-25}).  

Then according to lemma \ref{ref-3.2.9-20} 
 $f_\ast\Pscr$ is  a local Cohen-Macauly sheaf on $X$. Hence by lemma
 \ref{ref-4.2.1-52}  
 $f_\ast \Pscr$ corresponds to a reflexive $\Oscr_{Y^+}$-module
 $\Qscr^+$. We claim that $\Qscr^+$ is a local projective generator
 for ${}^0\Per(Y^+/X)$. By the characterization given in Proposition
 \ref{ref-3.2.7-18} it follows that it is sufficient to check
 this in the formal case, i.e. in the situation of \S\ref{ref-4.3-56}. 

So let us for a moment assume that we are in the formal
situation. Then according to Theorem \ref{ref-3.5.5-45} and
lemma \ref{ref-4.2.2-53}
$\Pscr=\oplus_i \Mscr_i^{\oplus a_i}$ with $a_i>0$ and $f_\ast
\Pscr=\oplus_{i} M_i^{\oplus a_i}$ and by Proposition \ref{ref-4.3.1-57}
this is equal to $\oplus_{i} (N^+_i)^{\oplus a_i}$ and using lemma
\ref{ref-4.2.2-53} for $f^+$ this corresponds to
$\oplus_{i}=(\Nscr^+_i)^{\oplus a_i}$. Using lemma \ref{ref-4.2.1-52} for
$f^+$ we obtain $\Qscr^+=(\Nscr^+_i)^{\oplus a_i}$ and using Theorem
\ref{ref-3.5.5-45} for $f^+$ we obtain that $\Qscr^+$ is a
projective generator for ${}^0\Per(Y^+/X)$.

Now we revert to the global case.  Put $\Ascr=f_\ast
\uEnd_Y(\Pscr)$. According to lemma \ref{ref-4.2.1-52} we also have $\Ascr=f_\ast
\uEnd_Y(\Qscr^+)$. 

According to Proposition
\ref{ref-3.3.1-24} we have now equivalences of categories
\begin{equation}
\label{ref-4.5-60}
D^b(\coh(Y))\xrightarrow{Rf_\ast \uRHom_Y(\Pscr,-)} D^b(\coh(\Ascr)) 
\xrightarrow{(f^+)^{-1}(-)\Lotimes_{(f^+)^{-1}(\Ascr)}\Qscr^+} D^b(\coh(Y^+))
\end{equation}
and these equivalences restrict to equivalences
\[
{}^{-1}\Per(Y/X)\r \coh(\Ascr)\r {}^0\Per(Y^+/X)
\]
This finishes the proof.
\end{proof}
\begin{remarks}
It is easy to see that the constructed equivalence is independent of
the choice of $\Pscr$. This is  expected since Bridgeland
constructs a canonical equivalence between $D^b(\coh(Y))$ and
$D^b(\coh(Y^+))$ \cite{Br1}. 
%
\end{remarks}
\appendix
\section{Non-commutative crepant resolutions?}
\label{ref-A-61}
Let $X=\Spec R$ where $(R,m)$ is a complete local $k$-algebra with
$R/m=k$.  Assume that $R$ has a three-dimensional isolated Gorenstein terminal
singularity. According to \cite[Cor. 3.12]{Reid} $R$ is a deformation of a Du Val
singularity. Let $t$ be the deformation parameter.

If $f:Y\r X$ is a crepant resolution of singularities of $X$ then
according to Corollary \ref{ref-3.2.11-22} there exists a Cohen-Macaulay
module $M$ over $R$ such that $A=\End_R(M)$ is homologically
homogeneous.  To simplify the discussion below let us call an
$R$-algebra which is finitely generated and torsion free as $R$-module
an \emph{$R$-order}.

Now assume that $X$ does not necessarily have a crepant resolution of
singularities. A natural question is: does there always exist a
homologically homogeneous $R$-order $A$? The answer to this question
is yes for trivial reasons. If we put $R'=R[t^{1/n}]$ and $X'=\Spec
R'$ then according to \cite{Brieskorn} for some $n$, $X'$ has a
crepant resolution of 
singularities. So there is a homologically homogeneous $R'$-order
$A'$. It now suffices to take $A=A'_{R}$. 

Hence to make the question meaningful we have to establish some
rules. The most obvious condition to impose is that $A$ should be
\emph{crepant} over $R$. By this we mean that
$\omega_A=\Hom_R(A,\omega_R)$ is equal to $A\otimes_R \omega_R$
where the inclusion $A\otimes_R
\omega_R\r\Hom_R(A,\omega_R)\cong\Hom_R(A,R)\otimes_R\omega_R$ is
obtained from the inclusion $A\r \Hom_R(A,R):a\mapsto \Tr(a-)$ where
$\Tr$ denotes the reduced trace map. So the reduced trace map must be
non-degenerate which is the same as saying that $A$ should be Azumaya
in codimension one \cite{Reiner}.  Let $\eta\in X$ be the generic
point. We also want $A$ and $R$ to be ``birational'' in some sense.
If we take birational to mean that $k(\eta)$ and $A_\eta$ should have
the same module category then by Morita theory it follows that
$A_\eta$ should be a full matrix ring over $k(\eta)$. By \cite[Prop
4.2]{AG1} it follows that $A=\End_R(M)$ for a reflexive $R$-module $M$.

So the more restrictive question is: does there exist a reflexive
$R$-module $M$ such that $A=\End_R(M)$ is homologically homogeneous. 
Let us consider a concrete example.
\begin{example} 
\label{ref-A.1-62}
Let $ R=k[[x,y,z,t]]/(x^2+y^2+z^2-t^{2b+1})$. Then there is no
reflexive $R$-module $M$ such that $A=\End_R(M)$ is homologically homogeneous.
\end{example}
\begin{proof}
Assume that $M$ exists.  Let $x$ be the unique closed point in $X$. $X-\{x\}$
  is regular so it follows from \cite[Thm 4.3]{AG1} that 
  $M$ is projective in codimension two. 

We apply $\Hom_R(M,-)$ to the short exact sequence
\[
0\r M\xrightarrow{\times t} M \r M/tM \r 0
\]
This yields an exact sequence
\[
0\r A/t A \r \End_{R/tR}(M/tM)\r \Ext^1_R(M,M)
\]
$M$ is projective in codimension two so $\Ext^1_R(M,M)$  has finite length.

Since $M$ is reflexive, $M/tM$ is torsion free and furthermore the
cokernel of $(M/tM)\r (M/tM)^{\ast\ast}$ is finite dimensional. So we
have an inclusion
\[
\End_{R/tR}(M/tM)\hookrightarrow \End_{R/tR}((M/tM)^{\ast\ast})
\]
whose cokernel has finite length.

We conclude that the cokernel $A/tA\r
\End_{R/tR}((M/tM)^{\ast\ast})$  has finite length. Since $A$ is
Cohen-Macaulay, the same holds for $A/tA$. So $A/tA$ is reflexive and
hence $A/tA=\End_{R/tR}((M/tM)^{\ast\ast})$. 

Now write $\bar{M}=(M/tM)^{\ast\ast}$, $\bar{R}=R/tR$.
 Since $\bar{R}$ is
integrally closed of dimension two, the maximal Cohen-Macaulay modules
are precisely the reflexive modules. So $\bar{M}$ is Cohen-Macaulay.
The ring $\bar{R}$ has a simple singularity so its indecomposable maximal
Cohen-Macaulay modules are known. The only non-trivial one is the
ideal $\bar{I}=(z,x+iy)$. 
 So if $n=\rk M$ then $\bar{M}$ is the direct sum
of $n$ rank one Cohen-Macaulay modules over $\bar{R}$. The
corresponding idempotents in $A/tA$ may be lifted to
$A$. It follows that $M$ is also  the sum of $n$ rank one
reflexive modules. But since $R$ is factorial this implies that $M$ is
 free. 

So we conclude that $A\cong M_n(R)$. But then $A$
is Morita equivalent to $R$ and so it has infinite global
dimension. It follows that $A$ can not be homologically homogeneous.
\end{proof}
\begin{remark}
I think that for three-dimensional terminal
Gorenstein singularities the existence of commutative and
non-commutative crepant resolutions are equivalent. This can presumably be
shown with the same Fourier-Mukai method which was used to establish
the three-dimensional McKay correspondence. 
\end{remark}
\begin{remark}
If $ R=k[[x,y,z,t]]/(x^2+y^2+z^2-t^{2b+1})$ is as in example \ref{ref-A.1-62}
then there exist a homologically homogeneous $R$-order $A$ \emph{with
  center $R$}. However $A$ is not Azumaya in codimension one so it is
not crepant in the above sense.  This
example was independently discovered by Schofield. It can probably be
generalized to arbitrary Gorenstein terminal singularities.

Here is the construction of $A$: let
$S=k[[x,y,z,s]]/(x^2+y^2+z^2-s^{4b+2})$. Sending $s\mapsto -s$ defines
an automorphism $\sigma$ of order two of $S$ and if we put $t=s^2$
then $R=S^\sigma$. 

Let $I=(x+iy, z-s^{2b+1})\subset S$ and put $A'=\bigl(\begin{smallmatrix} S & I
  \\ I^{-1} & S\end{smallmatrix}\bigr)$. We have $\sigma(I)=(x+iy)
  I^{-1}$ and hence if $M=\bigl(\begin{smallmatrix} 0 & x+iy
  \\ 1 & 0\end{smallmatrix}\bigr)$ then $M\sigma(-)M^{-1}$ defines
  an automorphism $\tau$ of $A'$ of order two. We take $A=A'{*}G$
  where $G=\{1,\tau\}\cong \ZZ/2\ZZ$.  
\end{remark}
\begin{remark} If $R$ is as in lemma \ref{ref-A.1-62} then it is still
  possible that there exists some saturated (see \cite{Bondal4,BondalVdB})
  triangulated category $\Ascr$ which in some sense serves as a crepant
  resolution of $R$. Tom Bridgeland has shown me a heuristic argument
  which seems to indicate that such $\Ascr$ does not exist. It would be
  interesting to settle this matter.
\end{remark}

\begin{thebibliography}{10}

\bibitem{AV}
M.~Artin and J.-L. Verdier, {\em Reflexive modules over rational double
  points}, Math. Ann. {\bf 270} (1985), no.~1, 79--82.

\bibitem{AG1}
M.~Auslander and O.~Goldman, {\em Maximal orders}, Trans. Amer. Math. Soc. {\bf
  97} (1960), 1--24.

\bibitem{Bondal4}
A.~I. Bondal and M.~M. Kapranov, {\em Representable functors, {S}erre functors,
  and reconstructions}, Izv. Akad. Nauk SSSR Ser. Mat. {\bf 53} (1989), no.~6,
  1183--1205, 1337.

\bibitem{BO2}
A.~I. Bondal and D.~O. Orlov, {\em Derived categories of coherent
  sheaves (ICM talk)},
  available as math.AG/0206295.

\bibitem{Bondal1}
\bysame, {\em Semi-orthogonal decompositions for algebraic varieties},
  available as alg-geom/950601, 1996.

\bibitem{BondalVdB}
A.~I. Bondal and M.~Van~den Bergh, {\em Generators and representability of
  functors in commutative and non-commutative geometry}, available as
  math.AG/0204218.

\bibitem{Br1}
T.~Bridgeland, {\em Flops and derived categories}, Invent. Math. {\bf 147}
  (2002), 613--632.

\bibitem{BKR}
T.~Bridgeland, A.~King, and M.~Reid, {\em The {M}c{K}ay correspondence as an
  equivalence of derived categories}, J. Amer. Math. Soc. {\bf 14} (2001),
  no.~3, 535--554 (electronic).

\bibitem{Brieskorn}
E.~Brieskorn, {\em Die {A}ufl\"osung der rationalen {S}ingularit\"aten
  holomorpher {A}bbildungen}, Math. Ann. {\bf 178} (1968), 255--270.

\bibitem{BH}
K.~A. Brown and C.~R. Hajarnavis, {\em Homologically homogeneous rings}, Trans.
  Amer. Math. Soc. {\bf 281} (1984), 197--208.

\bibitem{CKM}
H.~Clemens, J.~Koll{\'a}r, and S.~Mori, {\em Higher-dimensional complex
  geometry}, Ast\'erisque (1988), no.~166, 144 pp. (1989).

\bibitem{Fossum}
R.~M. Fossum, {\em The divisor class group of a {K}rull domain},
  Springer-Verlag, New York, 1973, Ergebnisse der Mathematik und ihrer
  Grenzgebiete, Band 74.

\bibitem{EGAIII}
A.~Grothendieck, {\em \'{E}l\'ements de g\'eom\'etrie alg\'ebrique. {I}{I}{I}.
  \'{E}tude cohomologique des faisceaux coh\'erents. {I}}, Inst. Hautes
  \'Etudes Sci. Publ. Math. (1961), no.~11, 167.

\bibitem{HRS}
D.~Happel, I.~Reiten, and S.~Smalo, {\em Tilting in abelian categories and
  quasitilted algebra}, Memoirs of the AMS, vol. 575, Amer. Math. Soc., 1996.

\bibitem{RD}
R.~Hartshorne, {\em Residues and duality}, Lecture notes in mathematics,
  vol.~20, Springer Verlag, Berlin, 1966.

\bibitem{H}
\bysame, {\em Algebraic geometry}, Springer-Verlag, 1977.

\bibitem{Chen}
C.~Jiun-Cheng, {\em Flops and equivalences of derived categories for threefolds
  with only {G}orenstein singularities}, available as math.AG/0202005.

\bibitem{KV}
M.~Kapranov and E.~Vasserot, {\em Kleinian singularities, derived categories
  and {H}all algebras}, Math. Ann. {\bf 316} (2000), no.~3, 565--576.

\bibitem{Kawamata}
Y.~Kawamata, {\em ${D}$-and ${K}$-equivalence}, available as AG/0205287.

\bibitem{Ke1}
D.~Keeler, {\em Ample filters of invertible sheaves}, available as
  math.AG/0108068.

\bibitem{Keller1}
B.~Keller, {\em Deriving {DG}-categories}, Ann. Sci. {\'E}cole Norm. Sup. (4)
  {\bf 27} (1994), 63--102.

\bibitem{Kollar}
J.~Koll{\'a}r, {\em Flops}, Nagoya Math. J. {\bf 113} (1989), 15--36.

\bibitem{KM}
J.~Koll{\'a}r and S.~Mori, {\em Birational geometry of algebraic varieties},
  Cambridge University Press, Cambridge, 1998, With the collaboration of C. H.
  Clemens and A. Corti, Translated from the 1998 Japanese original.

\bibitem{Neeman1}
A.~Neeman, {\em The {G}rothendieck duality theorem via {B}ousfield's techniques
  and {B}rown representability}, J. Amer. Math. Soc. {\bf 9} (1996), 205--236.

\bibitem{Reid}
M.~Reid, {\em Young person's guide to canonical singularities}, Algebraic
  geometry, Bowdoin, 1985 (Brunswick, Maine, 1985) (Providence, RI), Amer.
  Math. Soc., Providence, RI, 1987, pp.~345--414.

\bibitem{Reiner}
I.~Reiner, {\em Maximal orders}, Academic Press, New York, 1975.

\end{thebibliography}
\ifx\undefined\bysame
\newcommand{\bysame}{\leavevmode\hbox to3em{\hrulefill}\,}
\fi

\end{document}